\newcount\ite\ite=1\def\0{\global\ite=1\1}
\def\1{\item{\rm(\romannumeral\the\ite)}\advance\ite1\quad}
\def\phi{\varphi}

\documentclass[12pt,twoside]{article}
\usepackage{amssymb}
\usepackage{setspace}

\font\teneufm=eufm10 scaled \magstep1
\font\seveneufm=eufm7 scaled \magstep1
\font\fiveeufm=eufm5  scaled \magstep1
\newfam\eufmfam

\textfont\eufmfam=\teneufm
\scriptfont\eufmfam=\seveneufm
\scriptscriptfont\eufmfam=\fiveeufm
\def\frak#1{{\fam\eufmfam\relax#1}}

\newfam\msbfam
\font\tenmsb=msbm10 scaled \magstep1  \textfont\msbfam=\tenmsb
\font\sevenmsb=msbm7 scaled \magstep1 \scriptfont\msbfam=\sevenmsb
\font\fivemsb=msbm5 scaled \magstep1  \scriptscriptfont\msbfam=\fivemsb

\def\dd#1{\raise1.5pt\hbox{$\,\partial\!$}/\raise-2.5pt\hbox{$\!\partial#1\,$}}

\def\tilde{\widetilde}
\def\hat{\widehat}
\def\5#1{{\mathcal #1}}

\def\CC{{\mathbb C}}

\def\PP{{\mathbb P}}

\def\ra{\rightarrow}

\def\GL{\mathop{\rm GL}\nolimits}
\def\SL{\mathop{\rm SL}\nolimits}

\def\Ann{\mathop{\rm Ann}\nolimits}

 \def\HollowBoxx #1#2#3{{\dimen0=#1 \advance\dimen0 by -#2
       \dimen1=#1 \advance\dimen1 by #3
        \vrule height 0pt depth #3 width #2
       \hskip -#3
       \vrule height #1 depth #3 width #3}}
 \def\LeftContraction{\mathord{\kern1.45pt \HollowBoxx{6pt}{3.5pt}{.4pt}}\,}

 \def\HollowBox #1#2#3{{\dimen0=#1 \advance\dimen0 by -#3
       \dimen1=#1 \advance\dimen1 by #3
        \vrule height #1 depth #3 width #3
        \vrule height 0pt depth #3 width #2
        \hskip -#3}}
 \def\RightContraction{\mathord{\, \HollowBox{6pt}{3.1pt}{.4pt}} \kern1.6pt}

\def\qed{{\hfill $\Box$}}
\newtheorem{theorem}{THEOREM}[section]

\newtheorem{lemma}[theorem]{Lemma}

\newtheorem{remark}[theorem]{Remark}
\newtheorem{proposition}[theorem]{Proposition}
\newtheorem{conjecture}[theorem]{Conjecture}

\begin{document}

\begin{center}
{\Large \bf Extracting Invariants of\\
\vspace{0.2cm}
Isolated Hypersurface Singularities\\
\vspace{0.4cm}
from their Moduli Algebras}
\vspace{0.4cm}\\
\normalsize M. G. Eastwood and A. V. Isaev
\end{center}

\begin{quotation} 
{\small \sl \noindent We use classical invariant theory to construct invariants of complex graded Gorenstein algebras of finite vector space dimension. As a consequence, we obtain a way of extracting certain numerical invariants of quasi-homogeneous isolated hypersurface singularities from their moduli algebras, which extends an earlier result due to the first author. Furthermore, we conjecture that the invariants so constructed solve the biholomorphic equivalence problem in the homogeneous case. The conjecture is easily verified for binary quartics and ternary cubics. We show that it also holds for binary quintics and sextics. In the latter cases the proofs are much more involved. In particular, we provide a complete list of canonical forms of binary sextics, which is a result of independent interest.}
\end{quotation}

\thispagestyle{empty}

\pagestyle{myheadings}
\markboth{M. G. Eastwood and A. V. Isaev}{Invariants of Isolated Hypersurface Singularities}

\setcounter{section}{0}

\section{Introduction}\label{intro}
\setcounter{equation}{0}

Let ${\mathcal O}_m$ be the local algebra of holomorphic function germs at the origin in $\CC^m$ with $m\ge 2$. For every hypersurface germ ${\mathcal V}$ at the origin (considered with the corresponding reduced complex structure) denote by $I({\mathcal V})$ the ideal of elements of ${\mathcal O}_m$ vanishing on ${\mathcal V}$. Let $f$ be a generator of $I({\mathcal V})$, and consider the complex associative commutative algebra ${\mathcal A}({\mathcal V})$ defined as the quotient of ${\mathcal O}_m$ by the ideal in ${\mathcal O}_m$ generated by $f$ and all its first-order partial derivatives. The algebra ${\mathcal A}({\mathcal V})$, called the {\it moduli algebra}\, or {\it Tjurina algebra}\, of ${\mathcal V}$, is independent of the choice of $f$ as well as the coordinate system near the origin; the moduli algebras of biholomorphically equivalent hypersurface germs are isomorphic. Clearly, ${\mathcal A}({\mathcal V})$ is trivial if and only if ${\mathcal V}$ is non-singular. Furthermore, it is well-known that $0<\dim_{\CC}{\mathcal A}({\mathcal V})<\infty$ if and only if the germ ${\mathcal V}$ has an isolated singularity (see, e.g. Chapter 1 in \cite{GLS}).

By a theorem due to Mather and Yau (see \cite{MY}), two hypersurface germs ${\mathcal V}_1$, ${\mathcal V}_2$ in $\CC^m$ with isolated singularities are biholomorphically equivalent if their moduli algebras ${\mathcal A}({\mathcal V}_1)$, ${\mathcal A}({\mathcal V}_2)$ are isomorphic. Thus, given the dimension $m$, the moduli algebra ${\mathcal A}({\mathcal V})$ determines ${\mathcal V}$ up to biholomorphism. In particular, if $\dim_{\CC}{\mathcal A}({\mathcal V})=1$, then ${\mathcal V}$ is biholomorphic to the germ of the hypersurface $\{z_1^2+\dots+z_m^2=0\}$, and if $\dim_{\CC}{\mathcal A}({\mathcal V})=2$, then ${\mathcal V}$ is biholomorphic to the germ of the hypersurface $\{z_1^2+\dots+z_{m-1}^2+z_m^3=0\}$. The proof of the Mather-Yau theorem does not provide an explicit procedure for recovering the germ ${\mathcal V}$ from the algebra ${\mathcal A}({\mathcal V})$ in general, and finding a way for reconstructing ${\mathcal V}$ (or at least some invariants of ${\mathcal V}$) from ${\mathcal A}({\mathcal V})$ is an interesting open problem. Motivated by this problem, we explicitly extract certain numerical invariants of ${\mathcal V}$ from ${\mathcal A}({\mathcal V})$ with $\dim_{\CC}{\mathcal A}({\mathcal V})>2$ under the assumption that the singularity of ${\mathcal V}$ is quasi-homogeneous.

Let ${\mathcal V}$ be a hypersurface germ having an isolated singularity. The singularity of ${\mathcal V}$ is said to be {\it quasi-homogeneous}\, if for some (hence for any) generator $f$ of $I({\mathcal V})$ there exist positive integers $p_1,\dots,p_m,q$ such that, modulo a biholomorphic change of coordinates near the origin, $f$ is the germ of a polynomial $Q$ satisfying $Q(t^{p_1}z_1,\dots,t^{p_m}z_m)\equiv t^qQ(z_1,\dots,z_m)$ for all $t\in\CC$. The singularity of ${\mathcal V}$ is said to be {\it homogeneous}\, if one can choose $Q$ to be homogeneous (i.e.~$p_1=\dots=p_m=1$). In the latter case, the condition $\dim_{\CC}{\mathcal A}({\mathcal V})>2$ implies $\deg Q\ge 3$. By a theorem due to Saito (see \cite{Sa1}), the singularity of ${\mathcal V}$ is quasi-homogeneous if and only if $f$ lies in the Jacobian ideal ${\mathcal J}(f)$ in ${\mathcal O}_m$, which is the ideal generated by all first-order partial derivatives of $f$. Hence, for a quasi-homogeneous singularity, ${\mathcal A}({\mathcal V})$ coincides with the {\it Milnor algebra}\, ${\mathcal O}_m/{\mathcal J}(f)$ for any generator $f$ of $I({\mathcal V})$. Therefore, if the singularity of ${\mathcal V}$ is quasi-homogeneous, the algebra ${\mathcal A}({\mathcal V})$ is a complete intersection ring, which implies that ${\mathcal A}({\mathcal V})$ is Gorenstein (see \cite{B}). Recall that a local complex commutative associative algebra ${\mathcal A}$ with $1<\dim_{\CC}{\mathcal A}<\infty$ is a {\it Gorenstein ring}\, if and only if for the annihilator $\Ann({\mathfrak m}):=\{u\in{\mathfrak m}: u\cdot{\mathfrak m}=0\}$ of its maximal ideal ${\mathfrak m}$ one has $\dim_{\CC}\Ann({\mathfrak m})=1$ (see, e.g.~\cite{H}). The property that ${\mathcal A}({\mathcal V})$ is Gorenstein characterizes quasi-homogeneous singularities (see, e.g.~\cite{Ma}). Next, if the singularity of ${\mathcal V}$ is quasi-homogeneous, the algebra ${\mathcal A}({\mathcal V})$ is (non-negatively) graded. More precisely, one has ${\mathcal A}({\mathcal V})=\oplus_{j\ge0}{\mathcal L}_{j}$, where ${\mathcal L}_{j}$ are linear subspaces of ${\mathcal A}({\mathcal V})$, with ${\mathcal L}_{j}{\mathcal L}_{k}\subset {\mathcal L}_{j+k}$ and ${\mathcal L}_0\simeq\CC$. The existence of such a grading on ${\mathcal A}({\mathcal V})$ also characterizes quasi-homogeneous singularities (see \cite{XY}).

A criterion for two complex graded Gorenstein algebras of finite vector space dimension greater than 2 to be isomorphic was given in \cite{FIKK} (see also \cite{FK2}, \cite{I2} for results on algebras over arbitrary fields of characteristic zero). The criterion is stated in terms of certain polynomials that were first introduced in \cite{FK1}. Indeed, as explained in Section \ref{sect1} below, to every complex Gorenstein algebra ${\mathcal A}$ with $2<\dim_{\CC}{\mathcal A}<\infty$ one can associate polynomials of a special form on ${\mathfrak n}:={\mathfrak m}/\hspace{-0.1cm}\Ann({\mathfrak m})={\mathfrak m}/{\mathfrak m}^{\nu}$, called {\it nil-polynomials}, of degree $\nu$ with vanishing constant and linear terms, where $\nu\ge 2$ is the nil-index of ${\mathfrak m}$. In \cite{FIKK} it was shown that two complex graded Gorenstein algebras ${\mathcal A}_1$, ${\mathcal A}_2$ are isomorphic if and only if some (hence any) nil-polynomials $P_1$, $P_2$ arising from ${\mathcal A}_1$, ${\mathcal A}_2$, respectively, are {\it linearly equivalent up to scale}, that is, there exist $c\in\CC^*$ and a linear isomorphism $L:{\mathfrak n}_1\to{\mathfrak n}_2$ such that $cP_1=P_2\circ L$. The homogeneous component $P^{[s]}$ of degree $s$ of any nil-polynomial $P$ is in fact a polynomial on ${\mathfrak n}/{\mathfrak n}^{\nu+2-s}\simeq{\mathfrak m}/{\mathfrak m}^{\nu+2-s}$ for $s=2,\dots,\nu$. It then follows that for every $s$ and for any absolute classical invariant ${\mathcal I}$ of forms (i.e.~homogeneous polynomials) of degree $s$ on ${\mathfrak m}/{\mathfrak m}^{\nu+2-s}$, the value ${\mathcal I}(P^{[s]})$ is invariantly defined. Thus, application of absolute classical invariants to homogeneous components of nil-polynomials yields numerical invariants of complex graded Gorenstein algebras. This conclusion is summarized in Theorem \ref{theorem1} (see Section \ref{sect1}).

For a hypersurface germ ${\mathcal V}$ having a quasi-homogeneous singularity, Theorem \ref{theorem1} provides a way of extracting a collection of invariants of ${\mathcal V}$ from its moduli algebra ${\mathcal A}({\mathcal V})$ if $\dim_{\CC}{\mathcal A}({\mathcal V})>2$, and these invariants arise from absolute classical invariants of forms of degree $s$ on ${\mathfrak m}({\mathcal V})/{\mathfrak m}({\mathcal V})^{\nu+2-s}$ for $s=2,\dots,\nu$, where ${\mathfrak m}({\mathcal V})$ is the maximal ideal of ${\mathcal A}({\mathcal V})$. The idea of constructing invariants of quasi-homogeneous singularities from their moduli algebras using classical invariant theory was first proposed in \cite{Ea}, where a certain form on ${\mathfrak m}({\mathcal V})/{\mathfrak m}({\mathcal V})^2$ defined invariantly up to scale was introduced. In fact, this form coincides, up to scale, with the highest-order homogeneous component $P^{[\nu]}$ of any nil-polynomial $P$ arising from ${\mathcal A}({\mathcal V})$. In particular, Theorem 2.1 of \cite{Ea} is contained in our Theorem \ref{theorem1} for $s=\nu$. The invariants of ${\mathcal V}$ corresponding to $s=\nu$ will play an important role in our treatment of homogeneous singularities below.  

It is not clear whether the collection of invariants supplied by Theorem \ref{theorem1} is complete in the sense that the invariants so obtained solve the biholomorphic equivalence problem for hypersurface germs with quasi-homogeneous singularity. Regarding the question of completeness, we propose a conjecture, which states that in the case of homogeneous singularities the invariants provided by Theorem \ref{theorem1} do form a complete set and, moreover, just the invariants corresponding to $s=\nu$ are sufficient for solving the biholomorphic equivalence problem in the homogeneous case (see Conjecture \ref{conj1} in Section \ref{sect2}). This last statement does not hold for general quasi-homogeneous singularities (see Remark \ref{remqhom}).

Let ${\mathcal V}$ be a hypersurface germ having a homogeneous singularity. Then one has $\dim_{\CC}{\mathfrak m}({\mathcal V})/{\mathfrak m}({\mathcal V})^2=m$. In some coordinates $z=(z_1,\dots,z_m)$ near the origin the germ ${\mathcal V}$ is defined by a form $Q(z)$ of degree $n$ (recall that $n\ge 3$), where $Q$ is {\it minimal}, i.e.~the germ of $Q$ at the origin generates the ideal $I({\mathcal V})$. It then follows from results in \cite{Sa2} that $\nu=m(n-2)$. Thus, for a homogeneous singularity, every nil-polynomial $P$ has degree $m(n-2)$, and its highest-order homogeneous component $P^{[m(n-2)]}$ is a form in $m$ variables. We say that any of the mutually proportional forms of degree $m(n-2)$ arising in this way is {\it associated}\, to the form $Q$. 

It is easy to see that two hypersurface germs ${\mathcal V}_1$, ${\mathcal V}_2$ defined by minimal forms $Q_1$, $Q_2$, respectively, are biholomorphically equivalent if and only if the forms $Q_1$, $Q_2$ are {\it linearly equivalent}, that is, there exists $C\in\GL(m,\CC)$ for which $Q_1(Cz)\equiv Q_2(z)$. As shown in Proposition \ref{equivarbitr} in Section \ref{sect2}, for minimal forms $Q_1$, $Q_2$ such that the singularities of ${\mathcal V}_1$, ${\mathcal V}_2$ are isolated, the linear equivalence problem is solved by absolute classical invariants. Accordingly, Conjecture \ref{conj1} states that one can recover all absolute invariants of forms of degree $n$ in $m$ variables from absolute invariants of forms of degree\linebreak $m(n-2)$ in $m$ variables by evaluating the latter for associated forms. If correct, Conjecture \ref{conj1} would yield an explicit algorithm for extracting a complete set of invariants of homogeneous hypersurface singularities from their moduli algebras, which would complement the Mather-Yau theorem in this case.

As we will see in Section \ref{sect2}, Conjecture \ref{conj1} is easily verified for binary quartics ($m=2$, $n=4$) and ternary cubics ($m=3$, $n=3$). In Sections \ref{sect3} and \ref{sect4} (see Theorems \ref{theorem2} and \ref{theorem3}) we show that Conjecture \ref{conj1} also holds for binary quintics ($m=2$, $n=5$) and binary sextics ($m=2$, $n=6$). In the two latter cases the proofs are much harder and require rather lengthy computations, most of which have been performed using Maple. Our proof of Theorem \ref{theorem3} relies, in particular, on a complete list of canonical forms of binary sextics that we provide in Theorem \ref{normal6}.  It is well-known that a generic binary sextic is linearly equivalent to a sextic in Sylvester Canonical Form (see (\ref{sylcf})). In Theorem \ref{normal6} we exactly determine the sextics that are excluded by Sylvester's genericity assumptions, which is a result of independent interest.

{\bf Acknowledgements.} We would like to thank A. Gorinov for suggesting a proof of Proposition \ref{equivarbitr} below. This work is supported by the Australian Research Council.

\section{Nil-polynomials and classical invariants}\label{sect1}
\setcounter{equation}{0}

Let ${\mathcal A}$ be a complex Gorenstein algebra with $2<\dim_{\CC}{\mathcal A}<\infty$ and ${\mathfrak m}$ the maximal ideal of ${\mathcal A}$. Further, let $\exp_2: {\mathfrak m}\ra {\mathfrak m}$ be the map
$$
\displaystyle\exp_2(u):=\sum_{s=2}^{\infty}\frac{1}{s!}u^{s},\quad u\in{\mathfrak m}.
$$
By Nakayama's lemma, ${\mathfrak m}$ is a nilpotent algebra, and therefore the above sum is in fact finite, with the highest-order term corresponding to $s=\nu$, where $\nu\ge 2$ is the {\it nil-index}\, of ${\mathfrak m}$ (i.e.~the largest of all integers $\mu$ for which ${\frak m}^{\mu}\ne 0$). Using the map $\exp_2$, one can associate to the algebra ${\mathcal A}$ a collection of polynomials of a special form. Let $\Ann({\mathfrak m}):=\{u\in{\mathfrak m}: u\cdot{\mathfrak m}=0\}={\mathfrak m}^{\nu}$ be the annihilator of ${\mathfrak m}$, and $\Pi$ a hyperplane in ${\mathfrak m}$ complementary to $\Ann({\mathfrak m})$. A $\CC$-valued polynomial $P$ on $\Pi$ is called a {\it nil-polynomial}\, if there exists a linear form $\omega:{\mathfrak m}\ra\CC$ such that $\omega(\Ann({\mathfrak m}))=\CC$, $\ker\omega=\Pi$, and $P=\omega\circ\exp_2|_{\Pi}$. Observe that $\deg P=\nu$. 

If $P_1$, $P_2$ are two nil-polynomials arising from Gorenstein algebras ${\mathcal A}_1$, ${\mathcal A_2}$ and hyperplanes $\Pi_1\subset{\mathfrak m}_1$, $\Pi_2\subset{\mathfrak m}_2$, respectively, then $P_1$, $P_2$ are called {\it linearly equivalent up to scale}\, if there exist $c\in\CC^*$ and a linear isomorphism $L:\Pi_1\ra \Pi_2$ such that $cP_1=P_2\circ L$. It then follows from results of \cite{FIKK} (see also \cite{FK1}, \cite{FK2}, \cite{I2}) that the map $\varphi:{\mathfrak m}_1\ra{\mathfrak m_2}$, $\phi(u+v):=L(u)+c\,\tilde\omega_2^{-1}(\omega_1(v))$, is an algebra isomorphism, where $u\in\Pi_1$, $v\in\Ann({\mathfrak m}_1)$, $\tilde\omega_2:=\omega_2|_{\Ann({\mathfrak m}_2)}$, and $\omega_1$, $\omega_2$ are the linear forms corresponding to $P_1$, $P_2$, respectively. Thus, if $P_1$, $P_2$ are linearly equivalent up to scale, the algebras ${\mathcal A}_1$, ${\mathcal A_2}$ are isomorphic.

A nil-polynomial $P=\omega\circ\exp_2|_{\Pi}$ arising from a Gorenstein algebra ${\mathcal A}$ extends to the polynomial $\tilde P:=\omega\circ\exp_2$ on all of ${\mathfrak m}$, and one has $\tilde P(v)=0$, $\tilde P(u+v)=\tilde P(u)$ for all $u\in{\mathfrak m}$, $v\in\Ann({\mathfrak m})$. Therefore, $\tilde P$ gives rise to a polynomial $\hat P$ on the quotient ${\mathfrak m}/\hspace{-0.1cm}\Ann({\mathfrak m})={\mathfrak m}/{\mathfrak m}^{\nu}$. Further, let $P^{[s]}$ be the homogeneous component of order $s$ of $P$ and $\tilde P^{[s]}$ its extension to ${\mathfrak m}$, with $s=2,\dots,\nu$. One has $\tilde P^{[s]}(v)=0$, $\tilde P^{[s]}(u+v)=\tilde P^{[s]}(u)$ for all $u\in{\mathfrak m}$, $v\in{\mathfrak m}^{\nu+2-s}$. Thus, $\tilde P^{[s]}$ gives rise to a polynomial $\hat P^{[s]}$ on the quotient ${\mathfrak m}/{\mathfrak m}^{\nu+2-s}$. It then follows that if two nil-polynomials $P_1$, $P_2$ arising from Gorenstein algebras ${\mathcal A}_1$, ${\mathcal A_2}$, respectively, are linearly equivalent up to scale, there exist $c\in\CC^*$ and algebra isomorphisms $L^{[s]}: {\mathfrak m}_1/{\mathfrak m}_1^{\nu+2-s}\ra {\mathfrak m}_2/{\mathfrak m}_2^{\nu+2-s}$ such that $c\hat P_1^{[s]}=\hat P_2^{[s]}\circ L^{[s]}$, where $\nu$ is the nil-index of each of ${\mathcal A}_1$, ${\mathcal A_2}$ (note that the nil-indices of ${\mathcal A}_1$, ${\mathcal A_2}$ coincide since $\deg P_1=\deg P_2$). Observe also that for any nil-polynomials $P$ and $P'$ arising from the same algebra the corresponding highest-order components $\hat P^{[\nu]}$ and $\hat P^{'[\nu]}$ coincide up to scale.

Next, we say that the algebra ${\mathcal A}$ is {\it (non-negatively) graded}\, if ${\mathcal A}=\oplus_{j\ge0}{\mathcal L}_{j}$, where ${\mathcal L}_{j}$ are linear subspaces of ${\mathcal A}$, with ${\mathcal L}_{j}{\mathcal L}_{k}\subset {\mathcal L}_{j+k}$ and ${\mathcal L}_0\simeq\CC$. In this case ${\mathfrak m}=\oplus_{j>0}{\mathcal L}_j$ and $\Ann({\mathfrak m})={\mathcal L}_d$ for $d:=\max\{j:{\mathcal L}_{j}\ne 0\}$. It is shown in \cite{FIKK} (see also \cite{FK2}, \cite{I2}) that if ${\mathcal A}$ is graded, then all nil-polynomials arising from ${\mathcal A}$ are linearly equivalent up to scale. Thus, the following theorem holds.

\begin{theorem}\label{mainold} \sl Let ${\mathcal A}_1$, ${\mathcal A_2}$ be graded Gorenstein algebras with\linebreak $2<\dim_{\CC}{\mathcal A}_j<\infty$, $j=1,2$, and $P_1$, $P_2$ some nil-polynomials arising from ${\mathcal A}_1$, ${\mathcal A_2}$, respectively. Assume that ${\mathcal A}_1$, ${\mathcal A_2}$ are isomorphic. Then there exist $c\in\CC^*$ and algebra isomorphisms $L^{[s]}: {\mathfrak m}_1/{\mathfrak m}_1^{\nu+2-s}\ra {\mathfrak m}_2/{\mathfrak m}_2^{\nu+2-s}$ such that\linebreak $c\hat P_1^{[s]}=\hat P_2^{[s]}\circ L^{[s]}$.
\end{theorem}
Theorem \ref{mainold} allows one to utilize classical invariant theory for constructing certain numerical invariants of graded Gorenstein algebras. We will now recall the definitions of relative and absolute classical invariants (see, e.g.~\cite{O} for details). 

Let $W$ be a finite-dimensional complex vector space and ${\mathcal Q}_W^n$ the linear space of forms (i.e.~homogeneous polynomials) of a fixed degree $n\ge 2$ on $W$. Define an action of $\GL(W)$ on ${\mathcal Q}_W^n$ by the formula
$$
(C,Q)\mapsto Q_C,\quad  Q_C(w):=Q(C^{-1}w),\,\,\hbox{where $C\in\GL(W)$, $Q\in{\mathcal Q}_W^n$, $w\in W$.}
$$
Two forms are said to be {\it linearly equivalent}\, if they lie in the same orbit with respect to this action. More generally, for a subgroup $G\subset\GL(W)$ we say that two forms are $G$-{\it equivalent}\, if they lie in the same $G$-orbit. An {\it invariant}\, (or {\it relative classical invariant}) of forms of degree $n$ on $W$ is a polynomial $I:{\mathcal Q}_W^n\ra\CC$ such that for any $Q\in{\mathcal Q}_W^n$ and any $C\in\GL(W)$ one has $I(Q)=(\det C)^kI(Q_C)$, where $k$ is a non-negative integer called the {\it weight}\, of $I$. It follows that $I$ is in fact homogeneous of degree $k\dim_{\CC}W/n$.

Next, for any two invariants $I$ and $\tilde I$, with $\tilde I\not\equiv 0$, the ratio $I/\tilde I$ yields a rational function on ${\mathcal Q}_W^n$ that is defined, in particular, at the points where $\tilde I$ does not vanish. If $I$ and $\tilde I$ have equal weights, this function does not change under the action of $\GL(W)$, and we say that $I/\tilde I$ is an {\it absolute invariant}\, (or {\it absolute classical invariant}) of forms of degree $n$ on $W$. If one fixes coordinates $z_1,\dots,z_m$ in $W$, then any element $Q\in{\mathcal Q}^n_W$ is written as
$$
Q(z_1,\dots,z_m)=\sum_{i_1+\dots+ i_m=n}\left(\begin{array}{c} n\\ i_1,\dots,i_m\end{array}\right) a_{i_1,\dots,i_m}z_1^{i_1}\dots z_m^{i_m},
$$
where $a_{i_1,\dots,i_m}\in\CC$. In what follows we will introduce a number of absolute invariants that will be defined in terms of the coefficients $a_{i_1,\dots,i_m}$. Observe that for any absolute invariant ${\mathcal I}$ so defined its value ${\mathcal I}(Q)$ is in fact independent of the choice of coordinates in $W$. When working in coordinates, we always assume that $W=\CC^m$ and identify $\GL(W)$ with $\GL(m,\CC)$.

Now, Theorem \ref{mainold} yields the following result. 

\begin{theorem}\label{theorem1} \sl Let ${\mathcal A}$ be a graded Gorenstein algebra such that\linebreak $2<\dim_{\CC}{\mathcal A}<\infty$, and $P$ a nil-polynomial arising from ${\mathcal A}$. Further, for a fixed $s\in\{2,\dots,\nu\}$, let $W$ be a complex vector space isomorphic to ${\mathfrak m}/{\mathfrak m}^{\nu+2-s}$ by means of a linear map $\psi: W\ra {\mathfrak m}/{\mathfrak m}^{\nu+2-s}$, and ${\mathcal I}$ an absolute invariant of forms of degree $s$ on $W$. Then the value ${\mathcal I}(\psi^*\hat P^{[s]})$ depends only on ${\mathcal I}$ and the isomorphism class of ${\mathcal A}$.  
\end{theorem}

Let ${\mathcal V}$ be a hypersurface germ having a quasi-homogeneous singularity. Consider the moduli algebra ${\mathcal A}({\mathcal V})$ of ${\mathcal V}$, and let ${\mathfrak m}({\mathcal V})$ be the maximal ideal of ${\mathcal A}({\mathcal V})$. Assume that $\dim_{\CC}{\mathcal A}({\mathcal V})>2$. Theorem \ref{theorem1} allows one explicitly to extract from ${\mathcal A}({\mathcal V})$ invariants of ${\mathcal V}$ of the form ${\mathcal I}(\hat P^{[s]})$, where $P$ is a nil-polynomial arising from ${\mathcal A}({\mathcal V})$ and ${\mathcal I}$ is an absolute invariant of forms of degree $s$ on ${\mathfrak m}({\mathcal V})/{\mathfrak m}({\mathcal V})^{\nu+2-s}$, with $s=2,\dots,\nu$. For $s=\nu$ these invariants were first considered in article \cite{Ea} where a certain form of degree $\nu$ on ${\mathfrak m}({\mathcal V})/{\mathfrak m}({\mathcal V})^2$ was constructed in an invariant way (up to scale). It is straightforward to see that this form coincides, up to scale, with $\hat P^{[\nu]}$, where $P$ is any nil-polynomial arising from ${\mathcal A}({\mathcal V})$. Thus, the invariants introduced in \cite{Ea} are exactly the invariants supplied by Theorem \ref{theorem1} for $s=\nu$. 

It is natural to ask whether the invariants extracted from ${\mathcal A}({\mathcal V})$ as above form a complete set, i.e.~whether they solve the biholomorphic equivalence problem for hypersurface germs with quasi-homogeneous singularity. Below we will attempt to answer this question in the homogeneous case.

\section{The homogeneous case}\label{sect2}
\setcounter{equation}{0}

From now on we will only consider homogeneous singularities. We say that a non-zero form $Q\in{\mathcal Q}_{\CC^m}^n$ is {\it minimal}\, if the germ of $Q$ at the origin generates the ideal $I({\mathcal V}_Q)$, where ${\mathcal V}_Q$ is the germ of the hypersurface $\{Q=0\}$. If $Q$ is a binary form (i.e.~$m=2$), then it can be written as a product of non-zero linear factors. Each linear factor defines a point in $\CC\PP^1$ called a {\it root}\, of $Q$, and the multiplicities of the factors are referred to as the {\it root multiplicities}. The minimality of $Q$ then means that each root has multiplicity one, i.e.~$Q$ is {\it square-free}. Observe that for minimal forms $Q_1$, $Q_2$ the hypersurface germs ${\mathcal V}_{Q_{{}_1}}$, ${\mathcal V}_{Q_{{}_2}}$ are biholomorphically equivalent if and only if $Q_1$, $Q_2$ are linearly equivalent. We will now show that for minimal forms that define hypersurface germs with isolated singularity the linear equivalence problem is solved by absolute classical invariants.

For a form $Q\in{\mathcal Q}_{\CC^m}^n$, let $\Delta(Q)$ be its discriminant (see Chapter 13 in \cite{GKZ} for the definition\footnote{The formulae for the discriminant that we use below in the cases $m=2,3$ differ from the one given in \cite{GKZ} by scalar factors.}). The discriminant is a relative classical invariant of degree $m(n-1)^{m-1}$. Set
\begin{equation}
X^n_m:=\{Q\in{\mathcal Q}_{\CC^m}^n:\Delta(Q)\ne 0\}.\label{setx}
\end{equation}
Then $Q$ lies in $X^n_m$ if and only if $Q$ is minimal and the singularity of ${\mathcal V}_Q$ is isolated. If $Q$ is a binary form, then $Q\in X^n_2$ if and only if $Q$ is non-zero and square-free.

\begin{proposition}\label{equivarbitr} \sl For $n\ge 3$ the orbits of the $\GL(m,\CC)$-action on $X_m^n$ are separated by absolute classical invariants of the kind
\begin{equation}
{\mathcal I}=\frac{I}{\Delta^{p}},\label{formmsss}
\end{equation} 
where $p$ is a non-negative integer and $I$ is a relative classical invariant.
\end{proposition}
 
\noindent {\bf Proof:}\footnote{This proof was suggested to us by A. Gorinov.} The set $X^n_m$ has the structure of an affine algebraic variety, and with respect to this structure the action of the complex reductive group $G:=\GL(m,\CC)$ on $X^n_m$ is algebraic. The stabilizers of the $G$-action are finite (see \cite{OS}). Therefore, the quotient of $X^n_m$ by this action coincides with the Hilbert quotient $Z:=X^n_m/\hspace{-0.1cm}/G$ (see, e.g.~\cite{K}, \cite{MFK}). On the quotient $Z$ one can introduce the structure of an affine algebraic variety in such a way that the quotient map $\pi: X^n_m\ra Z$ is an algebraic morphism  and $\pi^*: \CC[Z]\ra \CC[X^n_m]^{G}$ is an isomorphism, where $\CC[Z]$ is the algebra of regular functions on $Z$ and $\CC[X^n_m]^{G}$ is the algebra of $G$-invariant regular functions on $X^n_m$. Since the points of $Z$ are separated by elements of $\CC[Z]$, the $G$-orbits in $X^n_m$ are separated by elements of $\CC[X^n_m]^{G}$. When $X^n_m$ is embedded into ${\mathcal Q}_{\CC^m}^n$ as in (\ref{setx}), every element of $\CC[X^n_m]^{G}$ becomes the restriction to $X_m^n$ of an absolute invariant of the form (\ref{formmsss}). The proof is complete. \qed
\vspace{0.3cm}

\noindent In what follows the algebra of the restrictions to $X_m^n$ of absolute invariants of the form (\ref{formmsss}) is denoted by ${\mathcal I}_m^n$. By the Hilbert Basis Theorem, this algebra is finitely generated.

Let now $Q\in X_m^n$ with $n\ge 3$. For the moduli algebra ${\mathcal A}({\mathcal V}_Q)$ of the germ ${\mathcal V}_Q$ it then follows from results in \cite{Sa2} that $\nu=m(n-2)$. Furthermore, one has $\dim_{\CC}{\mathfrak m}({\mathcal V}_Q)/{\mathfrak m}({\mathcal V}_Q)^2=m$. Thus, every nil-polynomial $P$ arising from ${\mathcal A}({\mathcal V}_Q)$ has degree $m(n-2)$, and the corresponding polynomial $\hat P^{[m(n-2)]}$ is a form on an $m$-dimensional space. For any two nil-polynomials $P$, $P'$ the forms $\hat P^{[m(n-2)]}$, $\hat P^{'[m(n-2)]}$ coincide up to scale, and we say that any of the mutually proportional forms of degree $m(n-2)$ arising in this way is {\it associated}\, to $Q$. Thus, every invariant of ${\mathcal V}_Q$ provided by Theorem \ref{theorem1} for $s=\nu$ is given as ${\mathcal I}({\bf Q})$, where ${\mathcal I}$ is an absolute classical invariant of forms of degree $m(n-2)$ on ${\mathfrak m}({\mathcal V}_Q)/{\mathfrak m}({\mathcal V}_Q)^2$ and ${\bf Q}$ is a form associated to $Q$.

For convenience, we will now make a canonical choice of variables in ${\mathfrak m}({\mathcal V}_Q)/{\mathfrak m}({\mathcal V}_Q)^2$. Consider the factorization maps $\pi_1:{\mathcal O}_m\ra {\mathcal O}_m/{\mathcal J}(Q)={\mathcal A}({\mathcal V}_Q)$ and $\pi_2:{\mathfrak m}({\mathcal V}_Q)\ra{\mathfrak m}({\mathcal V}_Q)/{\mathfrak m}({\mathcal V}_Q)^2$ (here and below, when speaking about ${\mathcal J}(Q)$, we identify $Q$ with its germ at the origin). Let $e_j$ be the image of the germ of the coordinate function $z_j$ under the composition $\pi_2\circ\pi_1$, $j=1,\dots,m$. Clearly, the vectors $e_j$ form a basis in ${\mathfrak m}({\mathcal V}_Q)/{\mathfrak m}({\mathcal V}_Q)^2$, and we denote by $w_1,\dots,w_m$ the coordinates with respect to this basis. For an absolute classical invariant ${\mathcal I}$ of forms of degree $m(n-2)$ in the variables $w_1,\dots,w_m$ it is easy to observe that ${\mathcal I}({\bf Q})$ is rational when regarded as a function of $Q$, with ${\bf Q}$ associated to $Q\in X_m^n$. 

Let ${\mathcal R}_m^n$ denote the collection of all invariant rational functions on $X_m^n$ obtained in this way. Further, let $\hat {\mathcal I}_m^n$ be the algebra of the restrictions to $X_m^n$ of {\it all}\, absolute invariants of forms of degree $n$ on $\CC^m$. Note that ${\mathcal R}_m^n$ lies in $\hat {\mathcal I}_m^n$ (see Proposition 1 in \cite{DC}). We propose the following:

\begin{conjecture}\label{conj1} \rm ${\mathcal R}_m^n=\hat {\mathcal I}_m^n$.
\end{conjecture}

\noindent Since every element of $\hat {\mathcal I}_m^n$ can be represented as a ratio of two elements of ${\mathcal I}_m^n$ (see Proposition 6.2 in \cite{Mu}), Conjecture \ref{conj1} is equivalent to the statement ${\mathcal I}_m^n\subset{\mathcal R}_m^n$.

If Conjecture \ref{conj1} were confirmed, it would provide a procedure for extracting a set of invariants of homogeneous  singularities from their moduli algebras that solves the biholomorphic equivalence problem for such singularities. This would be a step towards understanding how a germ ${\mathcal V}$ can be explicitly recovered from its moduli algebra ${\mathcal A}({\mathcal V})$ in general.

\begin{remark}\label{remqhom}\rm Observe that for general quasi-homogeneous singularities the invariants supplied by Theorem \ref{theorem1} for $s=\nu$ do not form a complete set. Indeed, consider the family of curves
$$
V_t:=\left\{(z_1,z_2)\in\CC^2:z_1^4+tz_1^2z_2^3+z_2^6=0\right\},\quad t\ne\pm 2,
$$
from Example 3.5 in \cite{FIKK} and let ${\mathcal V}_t$ be the germ of $V_t$ at the origin. Then, upon identification of ${\mathfrak m}({\mathcal V}_{t_{{}_1}})/{\mathfrak m}({\mathcal V}_{t_{{}_1}})^2$ and ${\mathfrak m}({\mathcal V}_{t_{{}_2}})/{\mathfrak m}({\mathcal V}_{t_{{}_2}})^2$, the highest-order terms of any nil-polynomials arising from ${\mathcal A}({\mathcal V}_{t_{{}_1}})$ and ${\mathcal A}({\mathcal V}_{t_{{}_2}})$ are linearly equivalent for all $t_1,t_2\ne 0$. However, as shown in \cite{FIKK}, ${\mathcal V}_{t_{{}_1}}$ and ${\mathcal V}_{t_{{}_2}}$ are biholomorphically equivalent if and only if $t_1=\pm t_2$.  
\end{remark}

Note that for binary quartics ($m=2$, $n=4$) and ternary cubics\linebreak ($m=3$, $n=3$) one has $m(n-2)=n$, that is, in these cases $\deg{\bf Q}=\deg Q$ for any form ${\bf Q}$ associated to $Q$, whereas in all other situations one has $\deg{\bf Q}>\deg Q$. In each of these two exceptional cases Conjecture \ref{conj1} states that every element of ${\mathcal I}_m^n$ can be recovered from some (possibly different) absolute invariant of forms of the same degree by applying them to associated forms. We will now show that this indeed holds.

Let $m=2$, $n=4$. Every non-zero square-free binary quartic is linearly equivalent to a binary quartic of the form
$$
q_t(z_1,z_2):=z_1^4+tz_1^2z_2^2+z_2^4,\quad t\ne\pm 2
$$ 
(see pp.~277--279 in \cite{El}). Any form associated to $q_t$ is again a binary quartic and is proportional to
$$
{\bf q}_t(w_1,w_2):=tw_1^4-12w_1^2w_2^2+tw_2^4
$$
(see \cite{Ea}). For $t\ne 0,\pm 6$ the quartic ${\bf q}_t$ is square-free, in which case the original quartic $q_t$ is associated to ${\bf q}_t$. 

The algebra of classical invariants of binary quartics is generated by certain invariants ${\rm I}_2$ and ${\rm I}_3$, where the subscripts indicate the degrees (see, e.g.~pp.~101--102 in \cite{El}). For a binary quartic of the form
$$
Q(z_1,z_2)=a_4z_1^4+6a_2z_1^2z_2^2+a_0z_2^4
$$
the values of the invariants ${\rm I}_2$ and ${\rm I}_3$ are computed as follows: 
\begin{equation}
\begin{array}{l}
{\rm I}_2(Q)=a_0a_4+3a_2^2,\quad {\rm I}_3(Q)=a_0a_2a_4-a_2^3,
\end{array}\label{invarquartics}
\end{equation}
and $\Delta(Q)={\rm I}_2(Q)^3-27\,{\rm I}_3(Q)^2$. Define an absolute invariant of binary quartics as
\begin{equation}
{\rm J}:=\frac{{\rm I}_2^3}{\Delta}.\label{absquartics}
\end{equation}
The restriction ${\rm J}|_{X_2^4}$ generates the algebra ${\mathcal I}_2^4$, and we have
$$
{\rm J}(q_t)=\frac{(t^2+12)^3}{108(t^2-4)^2}.
$$

Consider another (very similar) absolute invariant of binary quartics:
$$
{\rm K}:=\frac{{\rm I}_2^3}{27\,{\rm I}_3^2}.
$$
Then one obtains ${\rm K}({\bf q}_t)={\rm J}(q_t)$, and therefore ${\rm K}({\bf Q})={\rm J}(Q)$ for any $Q\in X_2^4$ and any ${\bf Q}$ associated to $Q$. Thus, the absolute invariant ${\rm K}$ evaluated for associated quartics yields a generator of ${\mathcal I}_2^4$, which agrees with Conjecture \ref{conj1}.

Let $m=3$, $n=3$. Every ternary cubic with non-zero discriminant is linearly equivalent to a ternary cubic of the form
$$
c_t(z_1,z_2,z_3):=z_1^3+z_2^3+z_3^3+tz_1z_2z_3,\quad t^3+27\ne 0
$$ 
(see p.~401 in \cite{W}). Any form associated to $c_t$ is again a ternary cubic and is proportional to
$$
{\bf c}_t(w_1,w_2,w_3):=tw_1^3+tw_2^3+tw_3^3-18w_1w_2w_3
$$
(see \cite{Ea}). For $t\ne 0$, $t^3-216\ne 0$ one has $\Delta({\bf c}_t)\ne 0$, in which case the original cubic $c_t$ is associated to ${\bf c}_t$.

The algebra of classical invariants of ternary cubics is generated by certain invariants ${\tt I}_4$ and ${\tt I}_6$, where, as before, the subscripts indicate the degrees (see pp.~381--389 in \cite{El}). For a ternary cubic of the form
$$
Q(z_1,z_2,z_3)=az_1^3+bz_2^3+cz_3^3+6dz_1z_2z_3
$$
the values of the invariants ${\tt I}_4$ and ${\tt I}_6$ are computed as follows: 
$$
\begin{array}{l}
{\tt I}_4(Q)=abcd-d^4 ,\quad {\tt I}_6(Q)=a^2b^2c^2-20abcd^3-8d^6 ,
\end{array}
$$
and $\Delta(Q)={\tt I}_6^2+64{\tt I}_4^3$. Define an absolute invariant of ternary cubics as
$$
{\tt J}:=\frac{{\tt I}_4^3}{\Delta}.
$$
The restriction ${\tt J}|_{X_3^3}$ generates the algebra ${\mathcal I}_3^3$, and we have
$$
{\tt J}(c_t)=-\frac{t^3(t^3-216)^3}{110592(t^3+27)^3}.
$$
Observe that ${\tt J}(c_t)=j(Z_t)/110592$, where $j(Z_t)$ is the value of the $j$-invariant for the elliptic curve $Z_t$ in $\CC\PP^2$ defined by the cubic $c_t$ (see, e.g.~\cite{I1}).

Consider another absolute invariant of ternary cubics:
$$
{\tt K}:=\frac{1}{4096\,{\tt J}}.
$$
Then one obtains ${\tt K}({\bf c}_t)={\tt J}(c_t)$, and therefore ${\tt K}({\bf Q})={\tt J}(Q)$ for any $Q\in X_3^3$ and any ${\bf Q}$ associated to $Q$. Thus, the absolute invariant ${\tt K}$ evaluated for associated cubics yields a generator of ${\mathcal I}_3^3$, which again agrees with Conjecture \ref{conj1}.

As we have seen, verification of Conjecture \ref{conj1} for binary quartics and ternary cubics is not hard. The first non-trivial case for testing Conjecture \ref{conj1} is that of binary quintics ($m=2$, $n=5$). We shall deal with this case in the next section.

\section{The case of binary quintics}\label{sect3}
\setcounter{equation}{0}

The main result of this section is the following:

\begin{theorem}\label{theorem2} \sl Conjecture {\rm \ref{conj1}} holds for binary quintics.
\end{theorem}

To prove Theorem \ref{theorem2}, we will need some facts from the invariant theory of binary quintics. First, let $Q\in{\mathcal Q}^n_{\CC^2}$ be a binary form of any degree $n\ge 2$ written as
$$
Q(z_1,z_2)=\sum_{i=0}^n\left(\begin{array}{l} n\\ i\end{array}\right) a_iz_1^iz_2^{n-i},
$$
where $a_i\in\CC$. The form $Q$ can be represented as a product of linear terms
$$
Q(z_1,z_2)=\prod_{\nu=1}^n(b_{\nu}z_1-c_{\nu}z_2),\label{representation}
$$
for some $b_{\nu},c_{\nu}\in\CC$. The discriminant of $Q$ is then given by
$$
\Delta(Q)=\frac{(-1)^{n(n-1)/2}}{n^n}\prod_{1\le\alpha<\beta\le n}(b_{\alpha}c_{\beta}-b_{\beta}c_{\alpha})^2
$$ 
(see pp.~97--101 in \cite{El}). The discriminant is a relative invariant of degree $2(n-1)$, which is non-zero if and only if $Q$ is non-zero and square-free. Furthermore, if $a_n\ne 0$, the discriminant $\Delta(Q)$ can be computed as
$$
\displaystyle\Delta(Q)=\frac{R(Q,\partial Q/\partial z_1)}{n^n a_n},
$$
where for two forms $P$ and $S$ we denote by $R(P,S)$ their resultant (see p.~36 in \cite{O}).

Next, define the $n$th {\it transvectant}\, as
$$
(Q,Q)^{(n)}:=(n!)^2\sum_{i=0}^n(-1)^i\left(\begin{array}{c} n \\ i\end{array}\right)a_ia_{n-i}
$$
(see Chapter 5 in \cite{O}). The transvectant $(Q,Q)^{(n)}$ is an invariant of degree 2. It is identically zero if $n$ is odd, thus for any odd $n$ we consider the invariant $\left(Q^2,Q^2\right)^{(2n)}$, which has degree 4. Observe that for the relative invariant ${\rm I}_2$ of binary quartics defined in (\ref{invarquartics}) one has ${\rm I}_2(Q)=(Q,Q)^{(4)}/1152$. 

We now introduce an absolute invariant of binary forms of degree $n$ as follows:
\begin{equation}
\displaystyle J(Q):=\left\{\begin{array}{ll}\displaystyle\frac{\Bigl[(Q,Q)^{(n)}\Bigr]^{n-1}}{\Delta(Q)} & \hbox{if $n$ is even,}\\
\vspace{-0.3cm}\\
\displaystyle\frac{\Bigl[\left(Q^2,Q^2\right)^{(2n)}\Bigr]^{(n-1)/2}}{\Delta(Q)} & \hbox{if $n$ is odd.}
\end{array}\right.\label{invariantj}
\end{equation}
Below, the invariant $J$ will be used for $n=5$ and $n=6$ (up to scale). Notice that for the absolute invariant ${\rm J}$ of binary quartics defined in (\ref{absquartics}) one has ${\rm J}=J/1152^3$.

We will now concentrate on the case of binary quintics. For $Q\in{\mathcal Q}_{\CC^2}^5$ define the {\it canonizant}\, of $Q$ as
$$
\hbox{Can}(Q):=\det\left(\begin{array}{lll}
a_5z_1+a_4z_2 & a_4z_1+a_3z_2 & a_3z_1+a_2z_2\\
a_4z_1+a_3z_2 & a_3z_1+a_2z_2 & a_2z_1+a_1z_2\\
a_3z_1+a_2z_2 & a_2z_1+a_1z_2 & a_1z_1+a_0z_2
\end{array}\right)
$$
(see p.~274 in \cite{El}). Clearly, $\hbox{Can}(Q)$ is a binary cubic, and we let $I_{12}(Q):=-27\,\Delta\left(\hbox{Can}(Q)\right)$ (see p.~307 in \cite{El}). It turns out that $I_{12}$ is a relative invariant of binary quintics of degree 12 (which explains the notation). In addition to $J$, we now introduce two more absolute invariants:
\begin{equation}
K(Q):=\frac{I_{12}(Q)^2}{\Delta(Q)^3},\quad L(Q):=\frac{\left(Q^2,Q^2\right)^{(10)}I_{12}(Q)}{\Delta(Q)^2}. \label{invariantk}
\end{equation}

It is well-known that  the algebra ${\mathcal I}_2^5$ is generated by the restrictions of $J$, $K$, $L$ to $X_2^5$. Indeed, let $I_{18}$ be the classical invariant of binary quintics of degree 18, as defined on p.~309 in \cite{El}. We have
\begin{equation}
16 I_{18}^2=I_4I_8^4+8I_8^3I_{12}-2I_4^2I_8^2I_{12}-72I_4I_8I_{12}^2-432I_{12}^3+I_4^3I_{12}^2 \label{i18}
\end{equation}
(see p.~313 in \cite{El}), where
$$
\begin{array}{l}
\displaystyle I_4(Q):=\frac{\left(Q^2,Q^2\right)^{(10)}}{7200000\cdot 10!},\\
\vspace{-0.1cm}\\
\displaystyle I_8(Q):=\frac{I_4(Q)^2-\Delta(Q)}{128}.
\end{array}
$$
The invariants $I_4$, $I_8$, $I_{12}$, $I_{18}$ generate the algebra of classical invariants of binary quintics (see \cite{Sy2}). Identity (\ref{i18}) then implies that for every absolute classical invariant ${\mathcal I}$ of binary quintics of the form (\ref{formmsss}) the numerator $I$ is a polynomial in $I_4$, $I_8$, $I_{12}$, and hence the restriction ${\mathcal I}|_{X_2^5}$ is a polynomial in $J|_{X_2^5}$, $K|_{X_2^5}$ and $L|_{X_2^5}$.

\begin{remark} \rm It is not hard to give examples showing that none of the three pairs of $J|_{X_2^5}$, $K|_{X_2^5}$, $L|_{X_2^5}$ generates ${\mathcal I}_2^5$. Indeed, for the quintics 
$$
\varphi_t:=z_1^5+tz_1^4z_2+z_2^5,\quad \hbox{with $256 t^5+3125\ne 0$},
$$
the invariants $K$ and $L$ vanish, but $J(\varphi_t)$ is a non-constant function of $t$. Next, for the quintics 
$$
\psi_t:=z_1^5+5tz_1^4z_2+5z_1z_2^4/t+z_2^5,\quad\hbox{with $t\ne 0$ and $t^5\ne 7\pm 4\sqrt{3}$}, 
$$
the invariants $J$ and $L$ vanish, but $K(\psi_t)$ is a non-constant function of $t$. Finally, for the quintics
$$
\rho_t:=z_1^5/t+z_2^5/(1-t)+(z_1+z_2)^5,\quad\hbox{with $t\ne 0,1, (1\pm i\sqrt{3})/2$}, 
$$
the value of $J$ is independent of $t$, and one can find $t_1$, $t_2$ such that $K(\rho_{t_1})=K(\rho_{t_2})$, but $L(\rho_{t_1})=-L(\rho_{t_2})\ne 0$.
\end{remark}

We will now present some canonical forms of square-free binary quintics. The forms given in the following proposition differ from the well-known canonical forms (such as Hammond's form shown on p.~305 in \cite{El}), but they will be particularly suitable for our purposes.

\begin{proposition}\label{normal} \sl Every non-zero square-free binary quintic is linearly equivalent to either a quintic of the form
\begin{equation}
f_{s,t}:=z_1^5+sz_1^4z_2+tz_1^3z_2^2+z_2^5,\quad s,t\in\CC,\label{familyfst}
\end{equation}
or a quintic of the form
\begin{equation}
f_{t}:=z_1^4z_2+tz_1^3z_2^2+z_2^5,\quad t\in\CC.\label{familyft}
\end{equation}
\end{proposition}

\noindent {\bf Proof:} As shown in \cite{El} (see pp.~268--271 and 305--306), a non-zero square-free binary quintic $Q$ is linearly equivalent to either
\begin{equation}
a_5z_1^5+5a_4z_1^4z_2+z_2^5,\quad \hbox{for some $a_4,a_5\in\CC$}\label{form1}
\end{equation}
(if $I_{12}(Q)=0$) or
\begin{equation}
az_1^5+bz_2^5+c(z_1+z_2)^5,\quad \hbox{for some $a,b,c\in\CC\setminus\{0\}$}\label{form2}
\end{equation}
(if $I_{12}(Q)\ne 0$). Observe that any square-free quintic of the form (\ref{form1}) is linearly equivalent to either $f_0$ (if $a_5=0$) or $f_{s,0}$ for some $s$ (if $a_5\ne 0$). Thus, in order to prove the proposition, it is sufficient to consider only quintics of the form (\ref{form2}).

\begin{lemma}\label{lemma} \sl Any square-free binary quintic of the form {\rm (\ref{form2})} is linearly equivalent to a quintic of the form
\begin{equation}
a_5'z_1^5+5a_4'z_1^4z_2+10a_3'z_1^3z_2^2+z_2^5,\quad a_3',a_4',a_5'\in\CC.\label{form3}
\end{equation}
\end{lemma}

\noindent {\bf Proof of Lemma \ref{lemma}:} Let us apply a general linear transformation
$$
z_1\mapsto \alpha z_1+\beta z_2,\quad z_2\mapsto \gamma z_1+\delta z_2,\quad \alpha,\beta,\gamma,\delta\in\CC
$$
to a quintic $Q$ of the form (\ref{form2}). To prove the lemma, we need to show that one can choose $\alpha,\beta,\gamma,\delta$ satisfying $\alpha\delta-\beta\gamma\ne 0$ such that
\begin{equation}
\begin{array}{l}
a\alpha^2\beta^3+b\gamma^2\delta^3+c(\alpha+\gamma)^2(\beta+\delta)^3=0,\\
\vspace{-0.1cm}\\
a\alpha\beta^4+b\gamma\delta^4+c(\alpha+\gamma)(\beta+\delta)^4=0.
\end{array}\label{cond}
\end{equation}

Set $\alpha=b\delta^3/a$, $\beta=1$, $\gamma=1$, and let $\delta$ be a root of the polynomial
$$
b\delta^3+c\left(\frac{b}{a}\delta^3+1\right)(\delta+1)^3.
$$
It is easy to see that $\alpha,\beta,\gamma,\delta$ so chosen satisfy identities (\ref{cond}). To see that $\alpha\delta-\beta\gamma\ne 0$, we assume the contrary and obtain
$$
a=b\delta^4,\quad c=-\frac{b\delta^4}{(1+\delta)^4}
$$ 
(observe that $\delta\ne -1$). It then follows that $\Delta(Q)=0$, which contradicts the assumption that $Q$ is square-free. Thus, $\alpha,\beta,\gamma,\delta$ define a non-degenerate linear transformation, and the proof of the lemma is complete. \qed
\vspace{0.3cm}

Proposition \ref{normal} follows from Lemma \ref{lemma}, since any square-free quintic of the form (\ref{form3}) is linearly equivalent to either $f_t$ for some $t$ (if $a_5'=0$) or $f_{s,t}$ for some $s,t$ (if $a_5'\ne 0$).\qed
\vspace{0.3cm}

Thus, in order to prove Theorem \ref{theorem2}, we need to consider two families of binary quintics: $f_{s,t}$ given by (\ref{familyfst}) and $f_t$ given by (\ref{familyft}), where the parameter values are such that $f_{s,t}$, $f_t$ are square-free. 

We start with the family $f_{s,t}$. First, we compute $J(f_{s,t})$, $K(f_{s,t})$, $L(f_{s,t})$. For the numerator and denominator in formula (\ref{invariantj}) we have, respectively,
\begin{equation}
(f_{s,t}^2,f_{s,t}^2)^{(10)}=57600\cdot 10!(125-3st^2)\label{numerator5}
\end{equation}
and
\begin{equation}
\displaystyle\Delta(f_{s,t})\!=\!\frac{1}{3125}\left(256 s^5-1600s^3t-27s^2t^4+2250st^2+108t^5+3125\right)\!.\label{denominator5}
\end{equation}
Formulae (\ref{invariantj}), (\ref{numerator5}), (\ref{denominator5}) yield
\begin{equation}
\begin{array}{l}
\displaystyle J(f_{s,t})=5(1440000\cdot 10!)^2(125-3st^2)^2/\\
\vspace{-0.3cm}\\
\hspace{1.7cm}\left(256s^5-1600s^3t-27s^2t^4+2250st^2+108t^5+3125\right).
\end{array}\label{jspecial}
\end{equation}
Further, by a straightforward computation we obtain
\begin{equation}
\begin{array}{l}
I_{12}(f_{s,t})=-\displaystyle\frac{1}{10^{10}}\left(19200 s^6t^2-160000 s^4t^3-1120 s^3t^6+\right.\\
\vspace{-0.3cm}\\
\hspace{1.95cm}
\left.440000 s^2t^4+3600 st^7+27 t^{10}-400000 t^5\right).
\end{array}\label{i12}
\end{equation}
Formulae (\ref{invariantk}), (\ref{numerator5}), (\ref{denominator5}), (\ref{i12}) imply
$$
\begin{array}{l}
\displaystyle K(f_{s,t})=\displaystyle\frac{1}{3125\cdot 4^{10}}\left(19200 s^6t^2-160000 s^4t^3-1120 s^3t^6+\right.\\
\vspace{-0.3cm}\\
\hspace{1.8cm}\left.440000s^2t^4+3600 st^7+27 t^{10}-400000 t^5\right)^2/\\
\vspace{-0.3cm}\\
\hspace{1.8cm}\left(256s^5-1600s^3t-27s^2t^4+2250st^2+108t^5+3125\right)^3,\\
\vspace{-0.1cm}\\
\displaystyle L(f_{s,t})=-\frac{225\cdot 10!}{4}(125-3st^2)\left(19200 s^6t^2-160000 s^4t^3-\right.\\
\vspace{-0.3cm}\\
\hspace{1.7cm}\left.1120 s^3t^6+440000 s^2t^4+3600 st^7+27 t^{10}-400000 t^5\right)/\\
\vspace{-0.3cm}\\
\hspace{1.7cm}\left(256 s^5-1600s^3t-27s^2t^4+2250st^2+108t^5+3125\right)^2.
\end{array}
$$

Next, any form associated to the quintic $f_{s,t}$ is proportional to the following binary sextic:
$$
\begin{array}{l}
{\bf f}_{s,t}(w_1,w_2):=(160s^3-300 st-27t^4)w_1^6+(-1200s^2+81st^3+\\
\vspace{-0.3cm}\\
\hspace{2.6cm}1125t)w_1^5w_2+(-270s^2t^2+3750s+675t^3)w_1^4w_2^2+\\
\vspace{-0.3cm}\\
\hspace{2.6cm}(480s^3t-1650st^2-6250)w_1^3w_2^3+(-480s^4+\\
\vspace{-0.3cm}\\
\hspace{2.6cm}2100s^2t-1125t^2)w_1^2w_2^4+(240s^3+27s^2t^3-825st-\\
\vspace{-0.3cm}\\
\hspace{2.6cm}108t^4)w_1w_2^5+(-6s^3t^2-50s^2+24st^3+125t)w_2^6.
\end{array}
$$
The algebra of classical invariants of binary sextics is generated by certain invariants of degrees 2, 4, 6, 10, 15, where the invariant of degree 10 can be taken to be the discriminant (see, e.g.~\cite{Sy2} and pp.~322--325 in \cite{El}). In this paper we use the generators utilized in \cite{Ea}. We write them in bracket form as follows (for the bracket notation see, e.g., Chapter 6 in \cite{O}):
\begin{equation}
\begin{array}{l}
{\mathsf I}_2:=[1,2]^6,\\
\vspace{-0.1cm}\\
{\mathsf I}_4:=[1,2]^4[1,3]^2[2,4]^2[3,4]^4,\\
\vspace{-0.1cm}\\
{\mathsf I}_6:=[1,2]^4[1,6]^2[2,3]^2[3,4]^4[4,5]^2[5,6]^4,\\
\vspace{-0.1cm}\\
{\mathsf I}_{10}:=[1,2]^2[1,3]^2[1,10]^2[2,3]^4[4,5]^2[4,6]^2[4,10]^2[5,6]^4[7,8]^2[7,9]^2\\
\vspace{-0.3cm}\\
\hspace{1.15cm}[7,10]^2[8,9]^4,\\
\vspace{-0.1cm}\\
{\mathsf I}_{15}:=[1,2]^2 [1,3]^2 [1,4] [2,3]^4 [4,5]^2[4,9][5,6]^2[6,7]^2[6,8]^2[7,8]^4[9,1]\\
\vspace{-0.3cm}\\
\hspace{1.15cm}[9,10]^4[10,11]^2[11,12]^4[12,13]^2[13,14]^2[13,15]^2[14,15]^4,
\end{array}\label{invarsextics}
\end{equation}
where, as before, the subscripts indicate the degrees (observe that ${\mathsf I}_2(Q)=(Q,Q)^{(6)}/(6!)^2$). 

Consider two absolute invariants of binary sextics:
\begin{equation}
{\mathsf J}:=\frac{3}{5}\frac{{\mathsf I}_2^2}{{\mathsf I}_2^2-2{\mathsf I}_4},\quad {\mathsf K}:=759375\frac{{\mathsf I}_{10}^2}{({\mathsf I}_2^2-2{\mathsf I}_4)^5}.\label{newjk}
\end{equation}
In \cite{Ea} the values ${\mathsf J}({\bf f}_{s,t})$ and ${\mathsf K}({\bf f}_{s,t})$ were computed as follows:
\begin{equation}
\begin{array}{l}
\displaystyle {\mathsf J}({\bf f}_{s,t})=(125-3st^2)^2/\\
\vspace{-0.3cm}\\
\hspace{1.6cm}(256s^5-1600s^3t-27s^2t^4+2250st^2+108t^5+3125),\\
\vspace{-0.1cm}\\
\displaystyle {\mathsf K}({\bf f}_{s,t})=F(s,t)^2/\\
\vspace{-0.3cm}\\
\hspace{1.65cm}\left(256s^5-1600s^3t-27s^2t^4+2250st^2+108t^5+3125\right)^3,
\end{array}\label{bfjk}
\end{equation}
where
$$
\begin{array}{l}
F(s,t):=163200s^6 t^2+14800000s^5-2100000s^4t^3+5400s^3t^6-\\
\vspace{-0.3cm}\\
\hspace{1.9cm}92500000s^3 t+7425000s^2 t^4-52650st^7+116250000st^2+\\
\vspace{-0.3cm}\\
\hspace{1.9cm}729 t^{10}-4556250 t^5+312500000.
\end{array}
$$
Now, from (\ref{numerator5}), (\ref{denominator5}), (\ref{i12}) we obtain
\begin{equation}
\begin{array}{l}
\displaystyle F(s,t)=-27\cdot 10^{10}\, I_{12}(f_{s,t})+\frac{115625}{4608\cdot 10!}\Delta(f_{s,t})(f_{s,t}^2,f_{s,t}^2)^{(10)}+\\
\vspace{-0.3cm}\\
\displaystyle \hspace{1.8cm}\frac{5}{2(19200\cdot 10!)^3}\Bigl((f_{s,t}^2,f_{s,t}^2)^{(10)}\Bigr)^3,
\end{array}\label{expressF}
\end{equation}
and formulae (\ref{invariantj}), (\ref{invariantk}), (\ref{denominator5}), (\ref{jspecial}), (\ref{bfjk}), (\ref{expressF}) yield
\begin{equation}
\begin{array}{l}
\displaystyle {\mathsf J}({\bf f}_{s,t})=\frac{1}{5(1440000\cdot 10!)^2}J(f_{s,t}),\\
\vspace{-0.1cm}\\
{\mathsf K}({\bf f}_{s,t})= 2^{20}\,3^6\,5^5\, K(f_{s,t})+c_1 L(f_{s,t})  + c_2 J(f_{s,t}) L(f_{s,t})+\\
\vspace{-0.3cm}\\
\hspace{1.7cm}c_3 J(f_{s,t})^3+c_4 J(f_{s,t})^2+c_5 J(f_{s,t})
\end{array}\label{connection1}
\end{equation}
for some $c_1,\dots,c_5\in\CC$.

Consider now a third absolute invariant of binary sextics:
\begin{equation}
{\mathsf L}:=675\frac{{\mathsf I}_2{\mathsf I}_{10}}{({\mathsf I}_2^2-2{\mathsf I}_4)^3}.\label{newl}
\end{equation} 
We have
\begin{equation}
\begin{array}{l}
\displaystyle {\mathsf L}({\bf f}_{s,t})=(125-3st^2)F(s,t)/\\
\vspace{-0.3cm}\\
\hspace{1.65cm}(256s^5-1600s^3t-27s^2t^4+2250st^2+108t^5+3125)^2.
\end{array}\label{bfl}
\end{equation}
Formulae (\ref{invariantj}), (\ref{invariantk}), (\ref{numerator5}), (\ref{denominator5}), (\ref{expressF}), (\ref{bfl}) yield
\begin{equation}
{\mathsf L}({\bf f}_{s,t})= -\frac{12}{25\cdot 10!}\,L(f_{s,t})+c_6 J(f_{s,t})^2+c_7 J(f_{s,t})\label{connection2}
\end{equation}
for some $c_6,c_7\in\CC$. 

We will now perform analogous calculations for the family $f_t$. First, we compute $J(f_t)$, $K(f_t)$, $L(f_t)$. For the numerator and denominator in formula (\ref{invariantj}) we have, respectively,
\begin{equation}
(f_t^2,f_t^2)^{(10)}=-172800\cdot 10!\,t^2\label{numerator51}
\end{equation}
and
\begin{equation}
\displaystyle\Delta(f_t)=\frac{256-27t^4}{3125}.\label{denominator51}
\end{equation}
Formulae (\ref{invariantj}), (\ref{numerator51}), (\ref{denominator51}) yield
\begin{equation}
\displaystyle J(f_t)=\frac{5(4320000\cdot 10!)^2\,t^4}{256-27t^4}.\label{jspecial1}
\end{equation}
Further, we obtain
\begin{equation}
I_{12}(f_t)=-\displaystyle\frac{1}{15625\cdot10^{10}}\left(421875\,t^{10}-175\cdot 10^5\,t^6+3\cdot 10^8\, t^2\right).\label{i121}
\end{equation}
Formulae (\ref{invariantk}), (\ref{numerator51}), (\ref{denominator51}), (\ref{i121}) imply
$$
\begin{array}{l}
\displaystyle K(f_t)=\displaystyle\frac{\left(421875\,t^{10}-175\cdot 10^5\,t^6+3\cdot 10^8\, t^2\right)^2}{4^{10}\cdot 5^{17}(256-27t^4)^3},\\
\vspace{-0.1cm}\\
\displaystyle L(f_t)=\frac{27\cdot 10!\,t^2\left(421875\,t^{10}-175\cdot 10^5\,t^6+3\cdot 10^8\, t^2\right)}{2500(256-27t^4)^2}.
\end{array}
$$

Next, any form associated to the quintic $f_t$ is proportional to the following binary sextic:
$$
\begin{array}{l}
{\bf f}_t(w_1,w_2):=(27 t^4-160)w_1^6-81t^3w_1^5w_2+270t^2w_1^4w_2^2-\\
\vspace{-0.1cm}\\
\hspace{2.4cm}480tw_1^3w_2^3+480w_1^2w_2^4-27t^3w_1w_2^5+6t^2w_2^6.
\end{array}
$$
We will now find ${\mathsf J}({\bf f}_t)$, ${\mathsf K}({\bf f}_t)$,  ${\mathsf L}({\bf f}_t)$ where ${\mathsf J}$, ${\mathsf K}$, ${\mathsf L}$ are the absolute invariants of binary sextics defined in (\ref{newjk}), (\ref{newl}). We obtain
\begin{equation}
\begin{array}{l}
\displaystyle {\mathsf J}({\bf f}_t)=\frac{9t^4}{256-27t^4},\\
\vspace{-0.1cm}\\
\displaystyle {\mathsf K}({\bf f}_t)=\frac{G(t)^2}{(256-27t^4)^3},\\
\vspace{-0.1cm}\\
\displaystyle {\mathsf L}({\bf f}_t)=-\frac{3t^2G(t)}{(256-27t^4)^2},
\end{array}\label{bfjk1}
\end{equation}
with $G(t):=729t^{10}+5400t^6+163200$. 
Now, from (\ref{numerator51}), (\ref{denominator51}), (\ref{i121}) we get
\begin{equation}
\begin{array}{l}
\displaystyle G(t)=-27\cdot 10^{10}\, I_{12}(f_t)+\frac{115625}{4608\cdot 10!}\Delta(f_t)(f_t^2,f_t^2)^{(10)}+\\
\vspace{-0.3cm}\\
\displaystyle \hspace{1.8cm}\frac{5}{2(19200\cdot 10!)^3}\Bigl((f_t^2,f_t^2)^{(10)}\Bigr)^3
\end{array}\label{expressF1}
\end{equation}
(cf. (\ref{expressF})), and formulae (\ref{invariantj}), (\ref{invariantk}), (\ref{denominator51}), (\ref{jspecial1}), (\ref{bfjk1}), (\ref{expressF1}) yield
\begin{equation}
\begin{array}{l}
\displaystyle {\mathsf J}({\bf f}_t)=\frac{1}{5(1440000\cdot 10!)^2}J(f_t),\\
\vspace{-0.1cm}\\
{\mathsf K}({\bf f}_t)= 2^{20}\,3^6\,5^5\, K(f_t)+c_1 L(f_t)  + c_2 J(f_t) L(f_t)+\\
\vspace{-0.3cm}\\
\hspace{1.5cm}c_3 J(f_t)^3+c_4 J(f_t)^2+c_5 J(f_t),\\
\vspace{-0.1cm}\\
\displaystyle {\mathsf L}({\bf f}_t)=-\frac{12}{25\cdot 10!}\,L(f_t)+c_6 J(f_t)^2+c_7 J(f_t),
\end{array}\label{connection11}
\end{equation}
where the constants coincide with the corresponding constants in formulae (\ref{connection1}), (\ref{connection2}).

Proposition \ref{normal} and formulae (\ref{connection1}), (\ref{connection2}), (\ref{connection11}) imply that for any non-zero square-free binary quintic $Q$ the following holds:
\begin{equation}
\begin{array}{l}
\displaystyle {\mathsf J}({\bf Q})=\frac{1}{5(1440000\cdot 10!)^2}J(Q),\\
\vspace{-0.1cm}\\
{\mathsf K}({\bf Q})=  2^{20}\,3^6\,5^5\, K(Q)+c_1 L(Q)+c_2 J(Q) L(Q)+\\
\vspace{-0.3cm}\\
\hspace{1.6cm}c_3 J(Q)^3+c_4 J(Q)^2+c_5 J(Q),\\
\vspace{-0.1cm}\\
\displaystyle {\mathsf L}({\bf Q})=-\frac{12}{25\cdot 10!}\,L(Q)+c_6 J(Q)^2+c_7 J(Q),
\end{array}\label{finalform}
\end{equation}
where ${\bf Q}$ is any binary sextic associated to $Q$. Since the restrictions of the invariants $J$, $K$, $L$ to $X_2^5$ generate the algebra ${\mathcal I}_2^5$, so do the restrictions to $X_2^5$ of the invariants that appear in the right-hand side of (\ref{finalform}). This completes the proof of Theorem \ref{theorem2}.\qed
\vspace{0.3cm}

As we have seen, verification of Conjecture \ref{conj1} for binary quintics is rather demanding computationally. In the next section we will show that the conjecture also holds for binary sextics ($m=2$, $n=6$), in which case even harder computations will be required.

\section{The case of binary sextics}\label{sect4}
\setcounter{equation}{0}

In this section we prove the following:

\begin{theorem}\label{theorem3} \sl Conjecture {\rm \ref{conj1}} holds for binary sextics.
\end{theorem}

First, we introduce eight absolute invariants of binary sextics. They are defined in terms of the discriminant and the invariants ${\mathsf I}_2$, ${\mathsf I}_4$, ${\mathsf I}_6$ (see (\ref{invarsextics})):
\begin{equation}
\begin{array}{lllll}
\displaystyle{\mathsf M}:=\frac{{\mathsf I}_2^5}{\Delta},&\displaystyle {\mathsf N}:=\frac{{\mathsf I}_4^5}{\Delta^2}, &\displaystyle {\mathsf P}:=\frac{{\mathsf I}_6^5}{\Delta^3},&\displaystyle {\mathsf R}:=\frac{{\mathsf I}_2\,{\mathsf I}_4^2}{\Delta},\\
\vspace{-0.1cm}\\
\displaystyle \displaystyle{\mathsf S}:=\frac{{\mathsf I}_2^3\,{\mathsf I}_4}{\Delta},&\displaystyle {\mathsf T}:=\frac{{\mathsf I}_4\,{\mathsf I}_6}{\Delta}, &\displaystyle {\mathsf U}:=\frac{{\mathsf I}_2^2\,{\mathsf I}_6}{\Delta}, &\displaystyle {\mathsf V}:=\frac{{\mathsf I}_2\,{\mathsf I}_6^3}{\Delta^2}
\end{array}\label{eight}
\end{equation}
(observe that ${\mathsf M}=J/(6!)^{10}$).

It is well-known that the algebra ${\mathcal I}_2^6$ is generated by the restrictions of the absolute invariants ${\mathsf M}$, ${\mathsf N}$, ${\mathsf P}$, ${\mathsf R}$, ${\mathsf S}$, ${\mathsf T}$, ${\mathsf U}$, ${\mathsf V}$ to $X_2^6$. Indeed, recall that the algebra of classical invariants of binary sextics is generated by the discriminant (which has degree 10) and ${\mathsf I}_2$, ${\mathsf I}_4$, ${\mathsf I}_6$, ${\mathsf I}_{15}$ (see (\ref{invarsextics})). Further, ${\mathsf I}_{15}^2$ is a polynomial in the discriminant and ${\mathsf I}_2$, ${\mathsf I}_4$, ${\mathsf I}_6$ (see \cite{Ea}). Therefore, for every absolute classical invariant ${\mathcal I}$ of binary quintics of the form (\ref{formmsss}) the numerator $I$ is a polynomial in the discriminant and ${\mathsf I}_2$, ${\mathsf I}_4$, ${\mathsf I}_6$. It is then straightforward to see that the restriction ${\mathcal I}|_{X_2^6}$ is a polynomial in the restrictions of the invariants ${\mathsf M}$, ${\mathsf N}$, ${\mathsf P}$, ${\mathsf R}$, ${\mathsf S}$, ${\mathsf T}$, ${\mathsf U}$, ${\mathsf V}$ to ${X_2^6}$.

Next, we present canonical forms of binary sextics. In order to prove Theorem \ref{theorem3}, it is sufficient to consider only the square-free case. However, the problem of finding a canonical form for {\it any}\, binary sextic is of independent interest, and we will give a solution to this problem in full generality. In \cite{Sy1} J. J. Sylvester showed that a generic binary sextic is linearly equivalent to a sextic of the following form (now called {\it Sylvester Canonical Form}):
\begin{equation}
au^6+bv^6+cw^6+duvw(u-v)(v-w)(w-u),\label{sylcf}
\end{equation}
for some $a,b,c,d\in\CC$ and linear forms $u$, $v$, $w$ satisfying $u+v+w=0$ any two of which are linearly independent. We take the opportunity here to explore Sylvester's proof with a view to determining exactly the sextics that are excluded by his genericity assumptions. Remarkably, it turns out that, up to linear equivalence, there are precisely seven sextics omitted from his canonical form.

\begin{theorem}\label{normal6} \sl A generic binary sextic is linearly equivalent to a sextic in Sylve\-ster Canonical Form. Any exceptional sextic is linearly equivalent to exactly one of the following seven sextics:
\begin{equation}
\makebox[250pt]{$\begin{array}{ll}
{\rm (i)} \,\,\,\,\,  z_1^4z_2^2, \quad {\rm (ii)}\,\,\, z_1^4(z_1^2+z_2^2), \quad {\rm (iii)} \,\,\, z_1^3z_2^3,\quad {\rm (iv)} \,\,\, z_1^5z_2,\\
\vspace{-0.1cm}\\
{\rm (v)} \,\,\,\, z_1(z_1^5+z_2^5),\quad {\rm (vi)} \,\,\,2z_1^6+18z_1^5z_2+10z_1^3z_2^3-z_2^6,\\
\vspace{-0.1cm}\\
{\rm (vii)} \, 184z_1^6-192z_1^5z_2-300z_1^4z_2^2-320z_1^3z_2^3-150z_1^2z_2^4-48z_1z_2^5+23z_2^6. 
\end{array}$}\label{listnormal6}
\end{equation}
\end{theorem}

\noindent {\bf Proof:} A general binary sextic
$$
Q=a_6z_1^6+6a_5z_1^5z_2+15a_4z_1^4z_2^2+20a_3z_1^3z_2^3+15a_2z_1^2z_2^4+6a_1z_1z_2^5+a_0z_2^6
$$
transforms under the $\SL(2,\CC)$-action on ${\mathcal Q}_{\CC^2}^6$ as the symmetric tensor $A_{ijklmn}$ given by
$$
\begin{array}{llll}
A_{111111}:=a_6, & A_{111112}:=a_5, & A_{111122}:=a_4, & A_{111222}:=a_3,\\
\vspace{-0.1cm}\\
A_{112222}:=a_2, & A_{122222}:=a_1, & A_{222222}:=a_0.
\end{array}
$$
Similarly, a general binary cubic
$$
Q'=a_3'z_1^3+3a_2'z_1^2z_2+3a_1'z_1z_2^2+a_0'z_2^3
$$
transforms under the $\SL(2,\CC)$-action on ${\mathcal Q}_{\CC^2}^3$ as the symmetric tensor $A'_{ijk}$ given by
$$
\begin{array}{llll}
A'_{111}:=a_3', & A'_{112}:=a_2', & A'_{122}:=a_1', & A'_{222}:=a_0'.
\end{array}
$$

For any $Q\in{\mathcal Q}_{\CC^2}^6$ we now consider the endomorphism $\hat Q: {\mathcal Q}_{\CC^2}^3\ra {\mathcal Q}_{\CC^2}^3$ defined by the following transformation of tensors:
$$
A'_{ijk}\mapsto
A_{ijklmn}\varepsilon^{lp}\varepsilon^{mq}\varepsilon^{nr}A'_{pqr},
$$
where $\varepsilon^{11}:=0$, $\varepsilon^{12}:=1$, $\varepsilon^{21}:=-1$, $\varepsilon^{22}:=0$. Observe that the map
\begin{equation}
{\mathcal Q}_{\CC^2}^6\times{\mathcal Q}_{\CC^2}^3\ra {\mathcal Q}_{\CC^2}^3,\quad (Q,Q')\mapsto \hat Q(Q')\label{mapPhi}
\end{equation}
is $\SL(2,\CC)$-equivariant. In terms of matrices, $\hat Q$ is given by
\begin{equation}\label{asmatrices}
\left(\begin{array}{c}z_1^3\\ z_1^2z_2\\ z_1z_2^2\\ z_2^3\end{array}\right)
\stackrel{\widehat Q}{\longmapsto}
\left(\begin{array}{rrrr}
-a_3 & -3a_2 & -3a_1 & -a_0\\
a_4 & 3a_3 & 3a_2 & a_1\\
-a_5 & -3a_4 & -3a_3& -a_2\\
a_6 & 3a_5 & 3a_4 & a_3\end{array}\right)
\left(\begin{array}{c}z_1^3\\ z_1^2z_2\\ z_1z_2^2\\ z_2^3\end{array}\right).
\end{equation}
The characteristic polynomial of the above matrix is the double quadratic
\begin{equation}
P(\lambda):=\lambda^4+\frac{{\mathsf I}_2(Q)}{2}\lambda^2-\frac{6{\mathsf I}_4(Q)-3{\mathsf I}_2(Q)^2}{16},\label{charpol}
\end{equation}
where ${\mathsf I}_2$ and ${\mathsf I}_4$ are the invariants defined in (\ref{invarsextics}). It may be solved explicitly to find 
the \lq eigencubics\rq\, of $\hat Q$, namely the binary cubics $Q'$ such that $\hat Q(Q')=\lambda Q'$ for some $\lambda\in{\mathbb{C}}$. 

Suppose first that $\hat Q$ has an eigencubic $Q'$ that is square-free. In this case $Q'$ can be chosen to be $\SL(2,\CC)$-equivalent to the cubic ${\mathsf c}:=z_1^3+z_2^3$. Owing to the $\SL(2,\CC)$-equivariance of map (\ref{mapPhi}), we then conclude that $Q$ is $\SL(2,\CC)$-equivalent to a sextic $Q_0$ satisfying $\hat Q_0({\mathsf c})=\lambda {\mathsf c}$. Now set
$$ 
Q_1:=Q_0-\lambda(z_1^6-z_2^6).
$$
Since $\widehat{z_1^6-z_2^6}({\mathsf c})={\mathsf c}$, it follows that $\hat Q_1({\mathsf c})=0$, which by formula (\ref{asmatrices}) occurs if and only if
$$
Q_1=a(z_1+z_2)^6+b(\omega z_1+\omega^2z_2)^6+c(\omega^2z_1+\omega z_2)^6
$$
for suitable $a,b,c\in\CC$, where $\omega:=e^{2\pi i/3}$. Setting
$$
u:=z_1+z_2,\quad v:=\omega z_1+\omega^2z_2,\quad w:=\omega^2z_1+\omega z_2
$$
and reassembling the result yields that $Q_0$ has the form (\ref{sylcf}). So far, we have merely recast Sylvester's original argument from \cite{Sy1} in modern language. His reasoning can be reversed, which shows that $Q$ is linearly equivalent to a sextic in Sylvester Canonical Form if and only if $\hat Q$ has a square-free eigencubic.

Now, we must deal with the case when none of the eigencubics of $\hat Q$ is square-free. There are two situations to consider: 
$$
\begin{array}{l}
\hbox{(a) $\hat Q$ has an eigencubic with a double root, and}\\
\vspace{-0.1cm}\\
\hbox{(b) $\hat Q$ has an eigencubic with a triple root.} 
\end{array}
$$
Passing to linearly equivalent forms, we may assume that the eigencubic is either $z_1^2z_2$ or $z_1^3$. In case it is $z_1^3$, we can read off from (\ref{asmatrices}) that $a_0=a_1=a_2=0$ and conclude that 
$$
\hat Q(z_1^3)=-a_3z_1^3\quad\mbox{and}\quad\hat Q(z_1^2z_2)=a_4z_1^3+3a_3z_1^2z_2.
$$
Therefore, we have
$$
\widehat Q(a_4z_1^3+4a_3z_1^2z_2)=3a_3(a_4z_1^3+4a_3z_1^2z_2).
$$
Hence if $a_3\ne 0$, then $a_4z_1^3+4a_3z_1^2z_2$ is an eigencubic with a double root. If $a_3=a_4=0$, then $z_1^2z_2$ is an eigencubic. Thus, in case (b) either we may reduce to case (a) or (if $a_3=0$, $a_4\ne 0$) by rescaling $z_2$ we may take $Q$ into a sextic of the form $a_6z_1^6+6\tilde a_5 z_1^5z_2+z_1^4z_2^2$. The substitution $z_2\mapsto z_2-3\tilde a_5z_1$ transforms this sextic into $\tilde a_6z_1^6+z_1^4z_2^2$. If $\tilde a_6=0$, we are led to $z_1^4z_2^2$. Otherwise, rescaling of $z_1$ and $z_2$ yields $z_1^4(z_1^2+z_2^2)$. We have thus obtained sextics (i) and (ii) in (\ref{listnormal6}). Observe that the eigencubics of each of $\widehat{z_1^4z_2^2}$, $z_1^2\widehat{(z_1^2+z_2^2)}$ are all proportional to $z_1^3$ (which has eigenvalue 0) and thus are not square-free.

From now on we assume that ${\mathsf c}_1:=z_1^2z_2$ is an eigencubic of $\hat Q$. Then (\ref{asmatrices}) implies $a_1=a_2=a_4=0$, and the characteristic polynomial $P(\lambda)$ given by formula (\ref{charpol}) factorizes as follows:
$$
P(\lambda)=(\lambda+3a_3)(\lambda-3a_3)(\lambda^2+a_0a_6-a_3^2).
$$
We may now explicitly find all the eigencubics. Clearly, ${\mathsf c}_1$ is an eigencubic with 
eigenvalue $3a_3$. Next, the following identity holds:
$$
(\widehat Q+3a_3\,{\mathrm{Id}})(8a_3^2a_5z_1^3-a_0a_5^2z_1^2z_2+2a_3(a_0a_6+8a_3^2)z_1z_2^2+2a_0a_3a_5z_2^3)=0.
$$
This identity implies that
$$
{\mathsf c}_2:=8a_3^2a_5z_1^3-a_0a_5^2z_1^2z_2+2a_3(a_0a_6+8a_3^2)z_1z_2^2+2a_0a_3a_5z_2^3
$$
is an eigencubic with eigenvalue $-3a_3$ provided ${\mathsf c}_2\not\equiv 0$. Further, for any $\lambda\in\CC$ the following identity takes place:
$$
\begin{array}{l}
(\widehat Q-\lambda\,{\mathrm{Id}})
((a_3-\lambda)(3a_3-\lambda)z_1^3-3a_0a_5z_1^2z_2+a_0(3a_3-\lambda)z_2^3)=\\
\vspace{-0.1cm}\\
\hspace{8cm}(3a_3-\lambda)(\lambda^2+a_0a_6-a_3^2)z_1^3.
\end{array}
$$
Consider the cubic
$$
{\mathfrak c}_{\lambda}:=(a_3-\lambda)(3a_3-\lambda)z_1^3-3a_0a_5z_1^2z_2+a_0(3a_3-\lambda)z_2^3,
$$
and let $\pm\lambda_0$ be the two square roots of $a_3^2-a_0a_6$. It then follows that ${\mathfrak c}_{\pm\lambda_0}$ is an eigencubic with eigenvalue $\pm\lambda_0$ provided ${\mathfrak c}_{\pm\lambda_0}\not\equiv 0$. It is straightforward to see that the cubics ${\mathsf c}_1$, ${\mathsf c}_2$, ${\mathfrak c}_{\lambda_0}$, ${\mathfrak c}_{-\lambda_0}$ are linearly independent if and only if
\begin{equation}
a_0a_3\lambda_0(\lambda_0-3a_3)(\lambda_0+3a_3)\ne 0.\label{cond1}
\end{equation}

Suppose first that condition (\ref{cond1}) holds. Then $Q$ is linearly equivalent to a sextic in Sylvester Canonical Form if and only if at least one of the eigencubics ${\mathsf c}_1$, ${\mathsf c}_2$, ${\mathfrak c}_{\lambda_0}$, ${\mathfrak c}_{-\lambda_0}$ is square-free. Therefore, we should determine the sextics for which all these eigencubics have multiple roots. The discriminants of the eigencubics (other than ${\mathsf c}_1$) are given by the following formulae:
$$
\begin{array}{l}
\Delta({\mathsf c}_2)=\displaystyle-\frac{4}{27}a_3a_5\Bigl(a_0^4a_3a_5^3a_6^2-128a_0^3a_3^3a_5^3a_6-32768a_3^{10}-\\
\vspace{-0.3cm}\\
\hspace{1.5cm}1536a_0^2a_3^6a_6^2+2a_0^4a_5^6-64a_0^3a_3^4a_6^3-2816a_0^2a_3^5a_5^3-12288a_0a_3^8a_6\Bigr),\\
\vspace{-0.1cm}\\
\Delta({\mathfrak c}_{\lambda_0})=-a_0^2(3a_3-\lambda_0)\Bigl(4a_0^2a_5^3-27a_3^5+81a_3^4\lambda_0-90a_3^3\lambda_0^2+\\
\vspace{-0.3cm}\\
\hspace{1.7cm}46a_3^2\lambda_0^3-11a_3\lambda_0^4+\lambda_0^5\Bigr),\\
\vspace{-0.1cm}\\
\Delta({\mathfrak c}_{-\lambda_0})=-a_0^2(3a_3+\lambda_0)\Bigl(4a_0^2a_5^3-27a_3^5-81a_3^4\lambda_0-90a_3^3\lambda_0^2-\\
\vspace{-0.3cm}\\
\hspace{1.9cm}46a_3^2\lambda_0^3-11a_3\lambda_0^4-\lambda_0^5\Bigr).
\end{array}
$$
Remarkably, the system of equations
$$
\Delta({\mathsf c}_2)=0,\,\,\Delta({\mathfrak c}_{\lambda_0})=0,\,\,
\Delta({\mathfrak c}_{-\lambda_0})=0,\,\,\lambda_0^2=a_3^2-a_0a_6
$$
may be solved explicitly, and we obtain 
$$
a_0=\frac{32a_5^6(115-41\sqrt{7})}{729a_6^5},\quad a_3=\frac{2a_5^3(1-2\sqrt{7})}{27a_6^2},\quad \lambda_0=\pm\frac{2a_5^3\sqrt{-11+4\sqrt{7}}}{3a_6^2},
$$
where $a_5,a_6\ne 0$. Up to linear equivalence, we may now assume that $a_5=3$, $a_6=1$ to obtain the sextic
\begin{equation}
z_1^6+18z_1^5z_2+40(1-2\sqrt{7})z_1^3z_2^3+32(115-41\sqrt{7})z_2^6.\label{sextic8prime}
\end{equation}
We claim that sextic (\ref{sextic8prime}) is linearly equivalent to sextic (vii) in (\ref{listnormal6}). Indeed, both these sextics are square-free and for either of them each absolute invariant introduced in (\ref{eight}) takes the same value, as shown below.
$$
\begin{array}{|c|c|c|c|c|c|c|c|}
\hline
{\mathsf M} & {\mathsf N} & {\mathsf P} & {\mathsf R} & {\mathsf S} & {\mathsf T} & {\mathsf U} & {\mathsf V}\\\hline
\vspace{-0.3cm}&&&&&&&\\
\displaystyle -\frac{7}{2^23^4} & \displaystyle \frac{5^5}{2^43^87^3} & \displaystyle -\frac{7^3}{2^63^{12}} & \displaystyle -\frac{5^2}{2^23^47} & \displaystyle -\frac{5}{2^23^4} & \displaystyle -\frac{5}{2^23^4} & \displaystyle -\frac{7}{2^23^4} & \displaystyle \frac{7^2}{2^43^8}\\
\vspace{-0.3cm}&&&&&&&\\ &&&&&&&\\ \hline
\end{array}
$$
Since the algebra ${\mathcal I}_2^6$ is generated by the restrictions of ${\mathsf M}$, ${\mathsf N}$, ${\mathsf P}$, ${\mathsf R}$, ${\mathsf S}$, ${\mathsf T}$, ${\mathsf U}$, ${\mathsf V}$ to $X_2^6$, it follows that the two sextics are indeed linearly equivalent. Thus if condition (\ref{cond1}) holds, we obtain sextic (vii) in (\ref{listnormal6}).    

We will now study the situations in which condition (\ref{cond1}) fails. Suppose that $\lambda_0=\pm3a_3$ and $a_3\ne 0$ (hence $a_0\ne 0$, $a_6\ne 0$). Then if $a_5=0$, the square-free cubic $2a_3\, z_1^3-a_0z_2^3$ is an eigencubic of $\hat Q$, thus we assume that $a_5\ne 0$. In this case any eigencubic of $\hat Q$ is proportional to either ${\mathsf c}_1$ (which has eigenvalue $3a_3$) or ${\mathsf c}_2$ (with has eigenvalue $-3a_3$). The condition $\Delta({\mathsf c}_2)=0$ yields    
$$
a_0=-\frac{32a_5^6}{729a_6^5},\quad a_3=\frac{2a_5^3}{27a_6^2}.
$$
Up to linear equivalence, we may assume that $a_5=3$, $a_6=2$ to obtain sextic (vi) in (\ref{listnormal6}).

Next, suppose that $\lambda_0=0$ and $a_3\ne 0$ (hence $a_0\ne 0$, $a_6\ne 0$). In this case any eigencubic is proportional to either ${\mathsf c}_1$ (which has eigenvalue $3a_3$) or ${\mathsf c}_2$ (which has eigenvalue $-3a_3$) or 
${\mathfrak c}_0=3a_0(a_6z_1^3-a_5z_1^2z_2+a_3z_2^3)$ (which has eigenvalue 0). If $a_5=0$, the cubic ${\mathfrak c}_0$ is square-free, and we assume that $a_5\ne 0$. It is now not hard to see that the system of equations
$$
\Delta({\mathsf c}_2)=0,\quad \Delta({\mathfrak c}_0)=0
$$
has no solutions.

Assume now that $\lambda_0=0$, $a_3=0$. Then $a_0a_6=0$. If $a_5=0$, $a_6\ne 0$ or $a_0\ne 0$, $a_5=0$ or $a_0=a_5=a_6=0$, the sextic $Q$ is in Sylvester Canonical Form. If $a_5\ne 0$, $a_6\ne 0$ or $a_0=0$, $a_5\ne 0$, $a_6=0$, then $Q$ is linearly equivalent to $z_1^5z_2$, which is sextic (iv) in (\ref{listnormal6}). Observe that any eigencubic of $\widehat{z_1^5z_2}$ is a linear combination of $z_1^3$ and $z_1^2z_2$ (which have eigenvalue 0) and thus is not square-free. If $a_0\ne 0$, $a_5\ne0$ then $Q$ is linearly equivalent to $z_1(z_1^5+z_2^5)$, which is sextic (v) in (\ref{listnormal6}). Observe that any eigencubic of $z_1\widehat{(z_1^5+z_2^5)}$ is proportional to $z_1z_2^2$ (which has eigenvalue 0) and thus is not square-free.

Next, let $\lambda_0\ne 0$, $a_3=0$ (hence $a_0\ne 0$, $a_6\ne 0$). If $a_5=0$, the sextic $Q$ is in Sylvester Canonical Form, and we assume that $a_5\ne 0$. In this case any eigencubic is proportional to either ${\mathsf c}_1$ (which has eigenvalue 0) or ${\mathfrak c}_{\lambda_0}$ (which has eigenvalue $\lambda_0$) or ${\mathfrak c}_{-\lambda_0}$ (which has eigenvalue $-\lambda_0$). It is now straightforward to observe that the system of equations
$$
\Delta({\mathfrak c}_{\lambda_0})=0,\quad \Delta({\mathfrak c}_{-\lambda_0})=0
$$
has no solutions.

Suppose finally that $a_0=0$, $a_3\ne 0$. Then if $a_5\ne 0$, the square-free cubic ${\mathsf c_2}=8a_3^2(a_5 z_1^3+2a_3z_1z_2^2)$ is an eigencubic of $\hat Q$. Hence we assume that $a_5=0$. In this case, if $a_6\ne 0$ the square-free cubic $a_6z_1^3+2a_3z_2^3$ is an eigencubic of $\hat Q$, and we further assume that $a_6=0$. Then $Q$ is linearly equivalent to $z_1^3z_2^3$, which is sextic (iii) in (\ref{listnormal6}). Observe that any eigencubic of $\widehat{z_1^3z_2^3}$ is proportional to either ${\mathsf c}_1$ (which has eigenvalue $3/20$) or $z_1z_2^2$ (which has eigenvalue $-3/20$) or $z_2^3$ (which has eigenvalue $1/20$) or $z_1^3$ (which has eigenvalue $-1/20$) and thus is not square-free.

It now remains to show that all the sextics in (\ref{listnormal6}) are pairwise linearly non-equivalent. Observe first that sextics (i)--(iv) have multiple roots with pairwise distinct patterns of root multiplicities. Further, each of sextics (v)--(vii) is square-free, and for them the absolute invariant ${\mathsf M}$ takes the values $0$, $9/637$, $-7/324$, respectively. Thus, no two sextics on list (\ref{listnormal6})  are linearly equivalent, and the proof of the theorem is complete.\qed

\begin{remark}\label{afterthought} \rm The exceptional sextics appearing in Theorem \ref{normal6} form a rather curious list. The roots of (i)--(v) are evidently nicely arranged (and sextic (v) has the appealing property that all of the invariants in (\ref{eight}) vanish). The roots of (vi) and (vii) are more mysterious. In both cases one can normalize three of them to be at 0, 1, $\infty$ to discover that the other three lie on the unit circle but that their arrangement is not especially symmetrical. In case (vi) one finds two bizarre angles of $60.35554...^\circ$ and $42.08850...^\circ$ and in case (vii) one finds $37.55551...^\circ$ and $28.00880...^\circ$ describing the arrangement of the remaining roots on the unit circle.   
\end{remark}

By Theorem \ref{normal6}, in order to prove Theorem \ref{theorem3} we need to consider non-zero square-free sextics of the form (\ref{sylcf}) as well as sextics (v)--(vii) in (\ref{listnormal6}). Any form associated to a non-zero square-free binary sextic is a binary octavic. Associated binary octavics generically can be found using the following Maple program:
$$
\begin{array}{l}
\hbox{\tt\small \verb-# Here we compute an associated octavic, we call it g:-}\\
\vspace{-0.5cm}\\
\hbox{\tt\small \verb-with(Groebner):-}\\
\vspace{-0.5cm}\\
\hbox{\tt\small \verb-F:=[f,diff(f,x),diff(f,y)]:-}\\
\vspace{-0.5cm}\\
\hbox{\tt\small \verb-G:=Basis(F,plex(x,y)):-}\\
\vspace{-0.5cm}\\
\hbox{\tt\small \verb-g:=numer(factor(coeff(NormalForm((X*x+Y*y)^8,G,plex(x,y)),y^8))):-}
\end{array}
$$
where {\tt f(x,y)} is the binary form for which an associated form is sought. The above program works if ${\tt y}^8\not\in{\mathcal J}({\tt f})$. Otherwise, one needs to take the coefficient at an appropriate monomial of degree 8 in ${\tt x,y}$ in the last line of this program. The resulting formulae for the associated octavics appear in the appendix at the end of this paper.

The algebra of classical invariants of binary octavics is generated by certain invariants of degrees 2, 3, 4, 5, 6, 7, 8, 9, 10 (see, e.g.~\cite{Sy2}). For the purposes of this paper we only require the invariants of degrees 2, 3, 4, 5 given in bracket form as follows:
$$
\begin{array}{l}
{\mathbf I}_2:=[1,2]^8,\\
\vspace{-0.1cm}\\
{\mathbf I}_3:=[1,2]^4[1,3]^4[2,3]^4,\\
\vspace{-0.1cm}\\
{\mathbf I}_4:=[1,2]^4[1,3]^4[2,4]^4[3,4]^4,\\
\vspace{-0.1cm}\\
{\mathbf I}_5:=[1,2]^4[2,3]^4[3,4]^4[4,5]^4[1,5]^4,
\end{array}
$$
where the subscripts indicate the degrees (note that ${\mathbf I}_2(Q):=(Q,Q)^{(8)}/(8!)^2$). Using these invariants, we will now introduce further invariants of degrees 4, 8, 10, 12 by the formulae
$$
\begin{array}{l}
\hat{\mathbf I}_4:=2{\mathbf I}_4-{\mathbf I}_2^2,\\
\vspace{-0.1cm}\\
{\mathbf I}_8:=297\,{\mathbf I}_2^4-1188\,{\mathbf I}_2^2\,{\mathbf I}_4+1188\,{\mathbf I}_4^2+1536\,{\mathbf I}_3\,{\mathbf I}_5-1280\,{\mathbf I}_2\,{\mathbf I}_3^2,\\
\vspace{-0.1cm}\\
{\mathbf I}_{10}:=-60\,{\mathbf I}_2\,{\mathbf I}_3\,{\mathbf I}_5+36\,{\mathbf I}_5^2+25\,{\mathbf I}_2^2\,{\mathbf I}_3^2,\\
\vspace{-0.1cm}\\
{\mathbf I}_{12}:=-512\,{\mathbf I}_2\,{\mathbf I}_{10}-8\,\hat{\mathbf I}_4\,{\mathbf I}_8+27\,\hat{\mathbf I}_4^3.
\end{array}
$$

Next, consider eight absolute invariants of binary octavics:
$$
\begin{array}{llll}
\displaystyle{\mathbf M}:=-\frac{3^8}{2^{11}5}\frac{\hat{\mathbf I}_4^5}{{\mathbf I}_{10}^2}, & \displaystyle{\mathbf N}:=\frac{3}{2^{22}5^{12}}\frac{{\mathbf I}_8^5}{{\mathbf I}_{10}^4}, & \displaystyle{\mathbf P}:=\frac{3^4}{2^{33}5^{18}}\frac{{\mathbf I}_{12}^5}{{\mathbf I}_{10}^6},&\displaystyle{\mathbf R}:=-\frac{3^2}{2^{11}5^5}\frac{\hat{\mathbf I}_4\,{\mathbf I}_8^2}{{\mathbf I}_{10}^2},\\
\vspace{-0.1cm}\\
\displaystyle{\mathbf S}:=-\frac{3^5}{2^{11}5^3}\frac{\hat{\mathbf I}_4^3\,{\mathbf I}_8}{{\mathbf I}_{10}^2}, & \displaystyle{\mathbf T}:=\frac{3}{2^{11}5^6}\frac{{\mathbf I}_8\,{\mathbf I}_{12}}{{\mathbf I}_{10}^2},& \displaystyle{\mathbf U}:=\frac{3^4}{2^{11}5^4}\frac{\hat{\mathbf I}_4^2\,{\mathbf I}_{12}}{{\mathbf I}_{10}^2}, & \displaystyle{\mathbf V}:=-\frac{3^4}{2^{22}5^{11}}\frac{\hat{\mathbf I}_4\,{\mathbf I}_{12}^3}{{\mathbf I}_{10}^4}.
\end{array}
$$
Let $Q$ be either a non-zero square-free binary sextic of the form (\ref{sylcf}) or one of sextics (v)--(vii) in (\ref{listnormal6}). Using Maple, one can now verify that for any binary octavic ${\mathbf Q}$ associated to $Q$ the following holds:
$$
\begin{array}{llll}
\displaystyle{\mathbf M}({\mathbf Q})={\mathsf M}(Q), & \displaystyle{\mathbf N}({\mathbf Q})={\mathsf N}(Q), & \displaystyle{\mathbf P}({\mathbf Q})={\mathsf P}(Q), & \displaystyle{\mathbf R}({\mathbf Q})={\mathsf R}(Q),\\
\vspace{-0.1cm}\\
\displaystyle{\mathbf S}({\mathbf Q})={\mathsf S}(Q), & \displaystyle{\mathbf T}({\mathbf Q})={\mathsf T}(Q), & \displaystyle{\mathbf U}({\mathbf Q})={\mathsf U}(Q), & \displaystyle{\mathbf V}({\mathbf Q})={\mathsf V}(Q).
\end{array}
$$
By Theorem \ref{normal6}, the above identities hold for any $Q\in X_2^6$. Since the algebra ${\mathcal I}_2^6$ is generated by the restrictions of ${\mathsf M}$, ${\mathsf N}$, ${\mathsf P}$, ${\mathsf R}$, ${\mathsf S}$, ${\mathsf T}$, ${\mathsf U}$, ${\mathsf V}$ to $X_2^6$, this completes the proof of Theorem \ref{theorem3}.\qed

\appendix

\setcounter{secnumdepth}{0}

\section{Appendix}

Here we give formulae for the octavics associated to non-zero square-free binary sextics of the form (\ref{sylcf}) as well as sextics (v)--(vii) in (\ref{listnormal6}). As stated in Section \ref{sect4}, these formulae were produced using a Maple program.

Let 
$$
Q=az_1^6+bz_2^6+c(z_1+z_2)^6+dz_1z_2(-z_1-z_2)(z_1-z_2)(z_1+2z_2)(-2z_1-z_2)
$$
be a sextics of the form (\ref{sylcf}). Then under the assumption $z_2^8\not\in{\mathcal J}(Q)$ any form associated to $Q$ is proportional to the octavic
$$
c_8w_1^8+c_7w_1^7w_2+c_6w_1^6w_2^2+c_5w_1^5w_2^3+c_4w_1^4w_2^4+c_3w_1^3w_2^5+c_2w_1^2w_2^6+c_1w_1w_2^7+c_0w_2^8,
$$
where the coefficients $c_j$ are given by the following expressions:
$$
\begin{array}{l}
c_0:=-2600c^3d^5-375c^4d^4+8110c^2d^6-9cd^7+9ad^7
-3150ac^4d^3-\\
\vspace{-0.1cm}\\
32970ac^3d^4+3656acd^6-8873ac^2d^5-378d^8
-175530a^2c^2d^4+8873a^2cd^5-\\
\vspace{-0.1cm}\\
8640a^2c^4d^2-95877a^2c^3d^3
+8110a^2d^6+297a^3b^2d^3-52596a^3c^3d^2-\\
\vspace{-0.1cm}\\
32970a^3cd^4-8991a^3c^4d
+95877a^3c^2d^3+2600a^3d^5+1750b^2cd^5+2429bcd^6+\\
\vspace{-0.1cm}\\
2601b^2c^2d^4+9887bc^2d^5-27b^3c^3d^2-297b^2c^3d^3-8070bc^3d^4
-540b^2c^4d^2+\\
\vspace{-0.1cm}\\
243ab^3c^4-81b^3c^4d+729a^2b^3c^3
+729a^3b^3c^2-27a^3b^3d^2-8991a^2b^2c^4+\\
\vspace{-0.1cm}\\
42768a^3b^2c^3
+2601a^2b^2d^4-1750ab^2d^5+17739a^3bc^4+1350bc^4d^3-\\
\vspace{-0.1cm}\\
9887a^2bd^5
+2429abd^6-81ab^3c^2d^2+162a^3b^3cd-81a^2b^3cd^2
-5103ab^2c^4d-\\
\vspace{-0.1cm}\\
25353a^2b^2c^3d-8370ab^2c^3d^2
+18090a^2b^2c^2d^2-8370a^3b^2cd^2+22797a^2b^2cd^3-\\
\vspace{-0.1cm}\\
13548ab^2cd^4
+26487a^2bc^4d+11070abc^4d^2+41661a^2bc^3d^2-375a^4d^4
-\\
\vspace{-0.1cm}\\
3402a^4c^4+250b^2d^6-13374abc^3d^3+41661a^3bc^2d^2
-73629abc^2d^4-\\
\end{array}
$$
$$
\hspace{-1.2cm}\begin{array}{l}
73629a^2bcd^4+8991a^4c^3d-8640a^4c^2d^2
+3150a^4cd^3+13374a^3bcd^3-\\
\vspace{-0.1cm}\\
162ab^3c^3d+25353a^3b^2c^2d
-22797ab^2c^2d^3+243a^4b^3c+81a^4b^3d-\\
\vspace{-0.1cm}\\
8991a^4b^2c^2+5103a^4b^2cd-540a^4b^2d^2+17739a^4bc^3-26487a^4bc^2d
+\\
\vspace{-0.1cm}\\
11070a^4bcd^2-1350a^4bd^3-8070a^3bd^4,
\end{array}
$$

$$
\begin{array}{l}
c_1:=8600c^3d^5+3000c^4d^4-36160c^2d^6-648cd^7-720ad^7
+23400ac^4d^3+\\
\vspace{-0.1cm}\\
151920ac^3d^4-14624acd^6+57944ac^2d^5
+1512d^8+702120a^2c^2d^4-\\
\vspace{-0.1cm}\\
13040a^2cd^5+57240a^2c^4d^2
+435456a^2c^3d^3-28720a^2d^6+1080a^3b^2d^3+\\
\vspace{-0.1cm}\\
210384a^3c^3d^2+111840a^3cd^4
+49248a^3c^4d-331560a^3c^2d^3-12200a^3d^5-\\
\vspace{-0.1cm}\\
14320b^2cd^5
-11096bcd^6+360b^3c^2d^3-45144b^2c^2d^4-125456bc^2d^5
+\\
\vspace{-0.1cm}\\
1728b^3c^3d^2+3456b^2c^3d^3+86520bc^3d^4+8640b^2c^4d^2
-5832ac^4b^3+648b^3c^4d-\\
\vspace{-0.1cm}\\
7776a^2b^3c^3+1944a^3b^3c^2
-1512a^3b^3d^2+360a^2b^3d^3+99144a^2b^2c^4-\\
\vspace{-0.1cm}\\
171072a^3b^2c^3
+24336a^2b^2d^4-320ab^2d^5-124416a^3bc^4-12600bc^4d^3-\\
\vspace{-0.1cm}\\
46360a^2bd^5
-8336abd^6-17496a^2b^3c^2d+1944ab^3c^2d^2-9072a^3b^3cd
-\\
\vspace{-0.1cm}\\
1296a^2b^3cd^2+720ab^3cd^3+66744ab^2c^4d
+269568a^2b^2c^3d+172800ab^2c^3d^2-\\
\vspace{-0.1cm}\\
72360a^2b^2c^2d^2
-105840a^3b^2cd^2+140976a^2b^2cd^3+54192ab^2cd^4-\\
\vspace{-0.1cm}\\
208656a^2bc^4d
-96120abc^4d^2+724464a^3bc^3d+521856a^2bc^3d^2+13608a^4c^4
-\\
\vspace{-0.1cm}\\
1000b^2d^6+3960bd^7+320112abc^3d^3-855144a^3bc^2d^2
+1651032a^2bc^2d^3+\\
\vspace{-0.1cm}\\
419256abc^2d^4+169776a^2bcd^4
-22680a^4c^3d+11880a^4c^2d^2-1800a^4cd^3+\\
\vspace{-0.1cm}\\
213120a^3bcd^3
-25056abcd^5-7776ab^3c^3d+66744a^3b^2c^2d+323352ab^2c^2d^3
+\\
\vspace{-0.1cm}\\
3888a^4b^3c-27216a^4b^2c^2+25920a^4b^2cd-4320a^4b^2d^2
-17496a^4bc^3+\\
\vspace{-0.1cm}\\
3240a^4bc^2d+7560a^4bcd^2-1800a^4bd^3-21960a^3bd^4,
\end{array}
$$

$$
\begin{array}{l}
c_2:=-11900c^3d^5-10500c^4d^4+85960c^2d^6+252cd^7+4536ad^7
-75600ac^4d^3-\\
\vspace{-0.1cm}\\
367920ac^3d^4+36512acd^6-194936ac^2d^5-3780d^8
-1282596a^2c^2d^4+\\
\vspace{-0.1cm}\\
37772a^2cd^5-162540a^2c^4d^2-936684a^2c^3d^3
+59920a^2d^6-61992a^3b^2d^3-\\
\vspace{-0.1cm}\\
343224a^3c^3d^2-227640a^3cd^4
-111132a^3c^4d+573048a^3c^2d^3+24500a^3d^5-\\
\end{array}
$$
$$
\begin{array}{l}
2100b^3cd^4
+66220b^2cd^5+32060bcd^6-13860b^3c^2d^3+120120b^2c^2d^4
+\\
\vspace{-0.1cm}\\
406364bc^2d^5-18900b^3c^3d^2+46116b^2c^3d^3-302820bc^3d^4
-49140b^2c^4d^2+\\
\vspace{-0.1cm}\\
47628ab^3c^4+2268b^3c^4d-61236a^2b^3c^3
-95256a^3b^3c^2-7560a^3b^3d^2+\\
\vspace{-0.1cm}\\
11340a^2b^3d^3-2100ab^3d^4
-435456a^2b^2c^4-598752a^3b^2c^3-123060a^2b^2d^4-\\
\vspace{-0.1cm}\\
14980ab^2d^5
+347004a^3bc^4+50400bc^4d^3+194992a^2bd^5+22400abd^6
+\\
\vspace{-0.1cm}\\
61236a^2b^3c^2d+68040ab^3c^2d^2-90720a^3b^3cd+79380a^2b^3cd^2
-2520ab^3cd^3-\\
\vspace{-0.1cm}\\
333396ab^2c^4d-1780380a^2b^2c^3d
-1068984ab^2c^3d^2-3283308a^2b^2c^2d^2-\\
\vspace{-0.1cm}\\
93744a^3b^2cd^2
-1118628a^2b^2cd^3-160440ab^2cd^4+687204a^2bc^4d
+\\
\vspace{-0.1cm}\\
355320abc^4d^2-2535624a^3bc^3d-2035908a^2bc^3d^2-20412a^4c^4
+9100b^2d^6-\\
\vspace{-0.1cm}\\
13860bd^7-928368abc^3d^3+2783592a^3bc^2d^2
-5778612a^2bc^2d^3-\\
\vspace{-0.1cm}\\
1237656abc^2d^4-364476a^2bcd^4+18144a^4c^3d
-3780a^4c^2d^2-937944a^3bcd^3+\\
\vspace{-0.1cm}\\
87696abcd^5+149688ab^3c^3d
+603288a^3b^2c^2d-506520ab^2c^2d^3+13608a^4b^3c
-\\
\vspace{-0.1cm}\\
4536a^4b^3d+6804a^4b^2c^2+9072a^4b^2cd-3780a^4b^2d^2
-27216a^4bc^3+\\
\vspace{-0.1cm}\\
31752a^4bc^2d-7560a^4bcd^2+76860a^3bd^4,
\end{array}
$$

$$
\begin{array}{l}
c_3:=8400c^3d^5+21000c^4d^4-140560c^2d^6-504cd^7-13104ad^7
+\\
\vspace{-0.1cm}\\
138600ac^4d^3+611520ac^3d^4-58352acd^6+448952ac^2d^5+6048d^8
+\\
\vspace{-0.1cm}\\
1390368a^2c^2d^4-
728a^2cd^5+257040a^2c^4d^2+1198008a^2c^3d^3
-70000a^2d^6-\\
\vspace{-0.1cm}\\
16632a^3b^2d^3+
293328a^3c^3d^2+252000a^3cd^4
+131544a^3c^4d-646632a^3c^2d^3-\\
\vspace{-0.1cm}\\
28000a^3d^5+
58800b^3cd^4
-81760b^2cd^5-56504bcd^6+136080b^3c^2d^3+\\
\vspace{-0.1cm}\\
79464b^2c^2d^4
-666008bc^2d^5+61992b^3c^3d^2-349272b^2c^3d^3+577920bc^3d^4
+\\
\vspace{-0.1cm}\\
143640b^2c^4d^2-
190512ab^3c^4-27216b^3c^4d+993384a^2b^3c^3
-503496a^3b^3c^2+\\
\vspace{-0.1cm}\\
16632a^3b^3d^2+
60480a^2b^3d^3-46200ab^3d^4
+993384a^2b^2c^4+2395008a^3b^2c^3+\\
\vspace{-0.1cm}\\
2184a^2b^2d^4+
112840ab^2d^5
-503496a^3bc^4-113400bc^4d^3-308728a^2bd^5-\\
\vspace{-0.1cm}\\
38864abd^6
+1832544a^2b^3c^2d-1144584ab^3c^2d^2-154224a^3b^3cd
+\\
\vspace{-0.1cm}\\
700056a^2b^3cd^2-433440ab^3cd^3+870912ab^2c^4d+2862216a^2b^2c^3d
+\\
\end{array}
$$
\vspace{0.3cm} 
$$
\begin{array}{l}
2549232ab^2c^3d^2+
10103184a^2b^2c^2d^2+704592a^3b^2cd^2
+897624a^2b^2cd^3+\\
\vspace{-0.1cm}\\
291648ab^2cd^4-1238328a^2bc^4d
-733320abc^4d^2+3510864a^3bc^3d+\\
\vspace{-0.1cm}\\
3234168a^2bc^3d^2
+
13608a^4c^4+7000b^3d^5-23800b^2d^6+22680bd^7+\\
\vspace{-0.1cm}\\
949536abc^3d^3
-4310712a^3bc^2d^2+8083152a^2bc^2d^3+2276232abc^2d^4
+\\
\vspace{-0.1cm}\\
468552a^2bcd^4-4536a^4c^3d+1352736a^3bcd^3-87696abcd^5
-879984ab^3c^3d-\\
\vspace{-0.1cm}\\
3578904a^3b^2c^2d-
1577016ab^2c^2d^3
+13608a^4b^3c-4536a^4b^3d+40824a^4b^2c^2-\\
\vspace{-0.1cm}\\
13608a^4b^2cd
+40824a^4bc^3-13608a^4bc^2d-126000a^3bd^4,
\end{array}
$$

$$
\begin{array}{l}
c_4:=-26250c^4d^4+165550c^2d^6+22050ad^7-157500ac^4d^3
-716100ac^3d^4+\\
\vspace{-0.1cm}\\
68348acd^6-676550ac^2d^5-7182d^8-788802a^2c^2d^4
-51100a^2cd^5-\\
\vspace{-0.1cm}\\
245700a^2c^4d^2-900900a^2c^3d^3+54250a^2d^6
+53550a^3b^2d^3-152712a^3c^3d^2-\\
\vspace{-0.1cm}\\
168000a^3cd^4-85050a^3c^4d
+431550a^3c^2d^3+17500a^3d^5-168000b^3cd^4+\\
\vspace{-0.1cm}\\
51100b^2cd^5
+68348bcd^6-431550b^3c^2d^3-788802b^2c^2d^4+676550bc^2d^5
-\\
\vspace{-0.1cm}\\
152712b^3c^3d^2+900900b^2c^3d^3-716100bc^3d^4-245700b^2c^4d^2
-13608a^3b^4c+\\
\vspace{-0.1cm}\\
411642ab^3c^4+85050b^3c^4d-2592324a^2b^3c^3
+1234926a^3b^3c^2-67662a^3b^3d^2-\\
\vspace{-0.1cm}\\
53550a^2b^3d^3+94500ab^3d^4
-1296162a^2b^2c^4-2592324a^3b^2c^3+12348a^2b^2d^4
-\\
\vspace{-0.1cm}\\
256900ab^2d^5+411642a^3bc^4+157500bc^4d^3+256900a^2bd^5
+48398abd^6-\\
\vspace{-0.1cm}\\
4734450a^2b^3c^2d+3406914ab^3c^2d^2-1233036a^2b^3cd^2
+1089900ab^3cd^3-\\
\vspace{-0.1cm}\\
1332450ab^2c^4d-3406536ab^2c^3d^2
-13513122a^2b^2c^2d^2-1233036a^3b^2cd^2-\\
\vspace{-0.1cm}\\
445704ab^2cd^4
+1332450a^2bc^4d+926100abc^4d^2-2438100a^3bc^3d
-\\
\vspace{-0.1cm}\\
3406536a^2bc^3d^2-3402a^4c^4-17500b^3d^5+54250b^2d^6-22050bd^7
-\\
\vspace{-0.1cm}\\
3402c^4b^4-13608ab^4c^3-20412a^2b^4c^2+3406914a^3bc^2d^2
-5761350a^2bc^2d^3-\\
\vspace{-0.1cm}\\
2781954abc^2d^4-445704a^2bcd^4
-1089900a^3bcd^3+2438100ab^3c^3d+\\
\vspace{-0.1cm}\\
4734450a^3b^2c^2d
+5761350ab^2c^2d^3-3402a^4b^4-13608a^4b^3c-20412a^4b^2c^2
-\\
\vspace{-0.1cm}\\
13608a^4bc^3+94500a^3bd^4,
\end{array}
$$

$$
\begin{array}{l}
c_5:=-8400c^3d^5+21000c^4d^4-140560c^2d^6+504cd^7-22680ad^7
+\\
\vspace{-0.1cm}\\
113400ac^4d^3+
577920ac^3d^4-56504acd^6+666008ac^2d^5+6048d^8
+\\
\vspace{-0.1cm}\\
79464a^2c^2d^4+81760a^2cd^5+
143640a^2c^4d^2+349272a^2c^3d^3
-23800a^2d^6-\\
\vspace{-0.1cm}\\
60480a^3b^2d^3+61992a^3c^3d^2+
58800a^3cd^4
+27216a^3c^4d-136080a^3c^2d^3-\\
\vspace{-0.1cm}\\
7000a^3d^5+252000b^3cd^4
+728b^2cd^5-58352bcd^6+646632b^3c^2d^3+\\
\vspace{-0.1cm}\\
1390368b^2c^2d^4
-448952bc^2d^5+293328b^3c^3d^2-1198008b^2c^3d^3+\\
\vspace{-0.1cm}\\
611520bc^3d^4
+257040b^2c^4d^2+13608a^3b^4c+4536a^3b^4d-
503496ab^3c^4
-\\
\vspace{-0.1cm}\\
131544b^3c^4d+2395008a^2b^3c^3-503496a^3b^3c^2+16632a^3b^3d^2
+
16632a^2b^3d^3-\\
\vspace{-0.1cm}\\
126000ab^3d^4+993384a^2b^2c^4+993384a^3b^2c^3
+2184a^2b^2d^4+
308728ab^2d^5-\\
\vspace{-0.1cm}\\
190512a^3bc^4-138600bc^4d^3
-112840a^2bd^5-38864abd^6+
13608ab^4c^2d+\\
\vspace{-0.1cm}\\
13608a^2b^4cd
+3578904a^2b^3c^2d-4310712ab^3c^2d^2+154224a^3b^3cd
+\\
\vspace{-0.1cm}\\
704592a^2b^3cd^2-1352736ab^3cd^3+1238328ab^2c^4d
-2862216a^2b^2c^3d+
\\
\vspace{-0.1cm}\\
3234168ab^2c^3d^2+10103184a^2b^2c^2d^2
+700056a^3b^2cd^2-897624a^2b^2cd^3+\\
\vspace{-0.1cm}\\
468552ab^2cd^4-870912a^2bc^4d
-733320abc^4d^2+879984a^3bc^3d+\\
\vspace{-0.1cm}\\
2549232a^2bc^3d^2
+28000b^3d^5-70000b^2d^6+13104bd^7+13608b^4c^4+\\
\vspace{-0.1cm}\\
40824ab^4c^3
+4536b^4c^3d+
40824a^2b^4c^2-949536abc^3d^3-1144584a^3bc^2d^2
+\\
\vspace{-0.1cm}\\
1577016a^2bc^2d^3+
2276232abc^2d^4+291648a^2bcd^4
+433440a^3bcd^3+\\
\vspace{-0.1cm}\\
87696abcd^5-3510864ab^3c^3d-
1832544a^3b^2c^2d
-8083152ab^2c^2d^3-46200a^3bd^4,
\end{array}
$$

$$
\begin{array}{l}
c_6:=-3780b^4c^2d^2+11900c^3d^5-10500c^4d^4+85960c^2d^6
-252cd^7+13860ad^7-\\
\vspace{-0.1cm}\\
50400ac^4d^3-302820ac^3d^4+32060acd^6
-406364ac^2d^5-3780d^8+\\
\vspace{-0.1cm}\\
120120a^2c^2d^4-66220a^2cd^5-49140a^2c^4d^2
-46116a^2c^3d^3+9100a^2d^6-\\
\vspace{-0.1cm}\\
11340a^3b^2d^3-18900a^3c^3d^2
-2100a^3cd^4-2268a^3c^4d+13860a^3c^2d^3-\\
\vspace{-0.1cm}\\
227640b^3cd^4
-37772b^2cd^5+36512bcd^6-573048b^3c^2d^3-1282596b^2c^2d^4
+\\
\vspace{-0.1cm}\\
194936bc^2d^5-343224b^3c^3d^2+936684b^2c^3d^3-367920bc^3d^4
-162540b^2c^4d^2+\\
\end{array}
$$
$$
\begin{array}{l}
13608a^3b^4c+4536a^3b^4d-3780a^2b^4d^2
+347004ab^3c^4+111132b^3c^4d-\\
\vspace{-0.1cm}\\
598752a^2b^3c^3-95256a^3b^3c^2
-7560a^3b^3d^2+61992a^2b^3d^3+76860ab^3d^4-\\
\vspace{-0.1cm}\\
435456a^2b^2c^4
-61236a^3b^2c^3-123060a^2b^2d^4-194992ab^2d^5+47628a^3bc^4
+\\
\vspace{-0.1cm}\\
75600bc^4d^3+14980a^2bd^5+22400abd^6-31752ab^4c^2d
-9072a^2b^4cd-\\
\vspace{-0.1cm}\\
7560ab^4cd^2-603288a^2b^3c^2d+2783592ab^3c^2d^2
+90720a^3b^3cd-93744a^2b^3cd^2+\\
\vspace{-0.1cm}\\
937944ab^3cd^3
-687204ab^2c^4d+1780380a^2b^2c^3d-2035908ab^2c^3d^2
-\\
\vspace{-0.1cm}\\
3283308a^2b^2c^2d^2+79380a^3b^2cd^2+1118628a^2b^2cd^3
-364476ab^2cd^4+\\
\vspace{-0.1cm}\\
333396a^2bc^4d+355320abc^4d^2-149688a^3bc^3d
-1068984a^2bc^3d^2-24500b^3d^5+\\
\vspace{-0.1cm}\\
59920b^2d^6-4536ad^7-20412b^4c^4
-27216ab^4c^3-18144b^4c^3d+6804a^2b^4c^2
+\\
\vspace{-0.1cm}\\
928368abc^3d^3+68040a^3bc^2d^2+506520a^2bc^2d^3-1237656abc^2d^4
-\\
\vspace{-0.1cm}\\
160440a^2bcd^4+2520a^3bcd^3-87696abcd^5+2535624ab^3c^3d
-61236a^3b^2c^2d+\\
\vspace{-0.1cm}\\
5778612ab^2c^2d^3-2100a^3bd^4,
\end{array}
$$

$$
\begin{array}{l}
c_7:=11880b^4c^2d^2+1800b^4cd^3-8600c^3d^5+3000c^4d^4
-36160c^2d^6+648cd^7-\\
\vspace{-0.1cm}\\
3960ad^7+12600ac^4d^3+86520ac^3d^4
-11096acd^6+125456ac^2d^5+1512d^8-\\
\vspace{-0.1cm}\\
45144a^2c^2d^4+14320a^2cd^5
+8640a^2c^4d^2-3456a^2c^3d^3-1000a^2d^6-360a^3b^2d^3+\\
\vspace{-0.1cm}\\
1728a^3c^3d^2
-648a^3c^4d-360a^3c^2d^3+111840b^3cd^4+13040b^2cd^5
-14624bcd^6+\\
\vspace{-0.1cm}\\
331560b^3c^2d^3+702120b^2c^2d^4-57944bc^2d^5
+210384b^3c^3d^2-435456b^2c^3d^3+\\
\vspace{-0.1cm}\\
151920bc^3d^4+57240b^2c^4d^2
+3888a^3b^4c-4320a^2b^4d^2+1800ab^4d^3-124416ab^3c^4
-\\
\vspace{-0.1cm}\\
49248b^3c^4d-171072a^2b^3c^3+1944a^3b^3c^2-1512a^3b^3d^2
-1080a^2b^3d^3-\\
\vspace{-0.1cm}\\
21960ab^3d^4+99144a^2b^2c^4-7776a^3b^2c^3
+24336a^2b^2d^4+46360ab^2d^5-5832a^3bc^4-\\
\vspace{-0.1cm}\\
23400bc^4d^3
+320a^2bd^5-8336abd^6-3240ab^4c^2d-25920a^2b^4cd+7560ab^4cd^2
-\\
\vspace{-0.1cm}\\
66744a^2b^3c^2d-855144ab^3c^2d^2+9072a^3b^3cd
-105840a^2b^3cd^2-213120ab^3cd^3+\\
\vspace{-0.1cm}\\
208656ab^2c^4d-269568a^2b^2c^3d
+521856ab^2c^3d^2-72360a^2b^2c^2d^2-1296a^3b^2cd^2
-\\
\vspace{-0.1cm}\\
140976a^2b^2cd^3+169776ab^2cd^4-66744a^2bc^4d
-96120abc^4d^2+7776a^3bc^3d+\\
\vspace{-0.1cm}\\
172800a^2bc^3d^2+12200b^3d^5
-28720b^2d^6+720bd^7+13608b^4c^4-17496ab^4c^3+\\
\end{array}
$$
$$
\begin{array}{l}
22680b^4c^3d
-27216a^2b^4c^2-320112abc^3d^3+1944a^3bc^2d^2
-323352a^2bc^2d^3+\\
\vspace{-0.1cm}\\
419256abc^2d^4+54192a^2bcd^4-720a^3bcd^3
+25056abcd^5-724464ab^3c^3d+\\
\vspace{-0.1cm}\\
17496a^3b^2c^2d-1651032ab^2c^2d^3,
\end{array}
$$

$$
\begin{array}{l}
c_8:=-375b^4d^4-8640b^4c^2d^2-3150b^4cd^3+2600c^3d^5
-375c^4d^4+8110c^2d^6+\\
\vspace{-0.1cm}\\
9cd^7-1350ac^4d^3-8070ac^3d^4+2429acd^6
-9887ac^2d^5-378d^8+2601a^2c^2d^4-\\
\vspace{-0.1cm}\\
1750a^2cd^5-540a^2c^4d^2
+297a^2c^3d^3+250a^2d^6-27a^3c^3d^2+81a^3c^4d-\\
\vspace{-0.1cm}\\
32970b^3cd^4
-8873b^2cd^5+3656bcd^6-95877b^3c^2d^3-175530b^2c^2d^4
+\\
\vspace{-0.1cm}\\
8873bc^2d^5-52596b^3c^3d^2+
95877b^2c^3d^3-32970bc^3d^4
-8640b^2c^4d^2+\\
\vspace{-0.1cm}\\
243a^3b^4c-81a^3b^4d-540a^2b^4d^2+1350ab^4d^3
+17739ab^3c^4+8991b^3c^4d+\\
\vspace{-0.1cm}\\
42768a^2b^3c^3+729a^3b^3c^2
-27a^3b^3d^2-297a^2b^3d^3-8070ab^3d^4-8991a^2b^2c^4
+\\
\vspace{-0.1cm}\\
729a^3b^2c^3+2601a^2b^2d^4+9887ab^2d^5+243a^3bc^4+3150bc^4d^3
+1750a^2bd^5+\\
\vspace{-0.1cm}\\
2429abd^6+26487ab^4c^2d-5103a^2b^4cd
+11070ab^4cd^2-25353a^2b^3c^2d+\\
\vspace{-0.1cm}\\
41661ab^3c^2d^2-162a^3b^3cd
-8370a^2b^3cd^2-13374ab^3cd^3-26487ab^2c^4d
+\\
\vspace{-0.1cm}\\
25353a^2b^2c^3d+41661ab^2c^3d^2+18090a^2b^2c^2d^2-81a^3b^2cd^2
-22797a^2b^2cd^3-\\
\vspace{-0.1cm}\\
73629ab^2cd^4+5103a^2bc^4d+11070abc^4d^2
+162a^3bc^3d-8370a^2bc^3d^2-\\
\vspace{-0.1cm}\\
2600b^3d^5+8110b^2d^6
-9bd^7-3402b^4c^4+17739ab^4c^3-8991b^4c^3d-8991a^2b^4c^2
+\\
\vspace{-0.1cm}\\
13374abc^3d^3-81a^3bc^2d^2+22797a^2bc^2d^3
-73629abc^2d^4-13548a^2bcd^4.
\end{array}
$$

Next, we consider sextics (v)--(vii) in (\ref{listnormal6}). For sextic (v) any associated form is proportional to
$$
28w_1^5w_2^3-3w_2^8,
$$
for sextic (vi) to
$$
\begin{array}{l}
35w_1^8+8w_1^7w_2-280w_1^6w_2^2-952w_1^5w_2^3+2072w_1^4w_2^4+\\
\vspace{-0.1cm}\\
\hspace{6cm}3080w_1^3w_2^5-2296w_1^2w_2^6
+944w_1w_2^7+224w_2^8,
\end{array}
$$
and for sextic (vii) to
$$
\begin{array}{l}
1148601562521w_1^8-3137415096128
w_1^7w_2+18448543636992w_1^6w_2^2+\\
\vspace{-0.1cm}\\
140944021374464w_1^5w_2^3-417766092163538w_1^4w_2^4+
282240199339072w_1^3w_2^5+\\
\vspace{-0.1cm}\\
74049821005648w_1^2w_2^6-25161150635776w_1w_2^7+18643008717596w_2^8.
\end{array}
$$

{\obeylines
\noindent Mathematical Sciences Institute
\noindent The Australian National University
\noindent Canberra, ACT 0200
\noindent Australia
\noindent e-mail: michael.eastwood@anu.edu.au, alexander.isaev@anu.edu.au
}

\end{document}